\documentclass[12pt]{amsart}
\usepackage{amsmath,amssymb}

\DeclareMathAlphabet{\curly}{U}{rsfs}{m}{n}

\newtheorem{thm}{Theorem}
\newtheorem{cor}{Corollary}
\newtheorem{lem}{Lemma}

\newtheorem{prop}{Proposition}

\theoremstyle{remark}
\newtheorem{remark}{Remark}

\newcommand{\RR}{{\mathbb R}}
\newcommand{\CC}{{\mathbb C}}
\newcommand{\NN}{{\mathbb N}}
\newcommand{\Fp}{{\mathbb F_p}}
\newcommand{\Fps}{{{\mathbb F}^*_p}}

\newcommand{\xx}{{\mathbf x}}
\newcommand{\yy}{{\mathbf y}}
\newcommand{\zz}{{\mathbf z}}
\newcommand{\uu}{{\mathbf u}}

\newcommand{\CA}{\curly A}
\newcommand{\CB}{\curly B}
\newcommand{\CP}{\curly P}

\newcommand{\cC}{\curly C}

\newcommand{\g}{\ensuremath{\gamma}}

\newcommand{\lam}{\ensuremath{\lambda}}
\newcommand{\eps}{\ensuremath{\varepsilon}}

\newcommand{\pfrac}[2]{{\left(\frac{#1}{#2}\right)}}

\newcommand{\be}{\begin{equation}}
\newcommand{\ee}{\end{equation}}
\newcommand{\benn}{\begin{equation*}}   
\newcommand{\eenn}{\end{equation*}}

\renewcommand{\b}{\ensuremath{\beta}}
\renewcommand{\a}{\ensuremath{\alpha}}

\renewcommand{\(}{\left(}
\renewcommand{\)}{\right)}

\newcommand{\order}{\asymp}

\newcommand{\CS}{Cauchy--Schwarz }
\newcommand{\cml}{\ensuremath{c_0}}  
\numberwithin{equation}{section}

\title[Explicit constructions of RIP matrices]
{Explicit constructions of RIP matrices and related problems}

\author[J. Bourgain]{Jean Bourgain}\address{School of Mathematics,
Institute for Advanced Study\\
Princeton, NJ 08540, USA} \email{bourgain@math.ias.edu}

\author{S. J. Dilworth}\address{Department of Mathematics \\
University of South Carolina\\ Columbia, SC 29208, USA}
\email{dilworth@math.sc.edu}

\author[K. Ford]{Kevin Ford} \address{Department of Mathematics\\ University
of Illinois at Urbana-Champaign\\1409 W. Green Street\\ Urbana, IL
61801, USA} \email{ford@math.uiuc.edu}

\author[S. Konyagin]{Sergei Konyagin}\address{Steklov Mathematical Institute\\
8, Gubkin Street, Moscow, 119991, Russia}
\email{konyagin@mi.ras.ru}

\author[D. Kutzarova]{Denka Kutzarova} \address{Institute of
Mathematics\\Bulgarian Academy of Sciences\\Sofia\\ Bulgaria}
\curraddr{Department of Mathematics\\
University of Illinois at Urbana-Champaign\\1409 W. Green Street\\
Urbana, IL 61801\\USA} \email{denka@math.uiuc.edu}
\date{\today}

\begin{document}
\thanks{The first author was partially supported by NSF Grants DMS 0808042 
and DMS 0835373, 
the second author was supported by NSF Grant DMS 0701552 and the
third author was supported by NSF grant DMS 0901339. 
The second, third and fifth authors were supported by the
Workshop in Analysis and Probability at Texas A\& M University, 2008.
The research
was finished while the third and fourth authors were visiting the
Institute for Advanced Study, the third author supported by grants
from the Ellentuck Fund and The Friends of the Institute For
Advanced Study. The last three authors thank the IAS for its
hospitality and excellent working conditions}

\begin{abstract}
We give a new explicit construction of $n\times N$ matrices
satisfying the Restricted Isometry Property (RIP).
Namely, for some $\eps>0$, large $N$ and any $n$ satisfying
$N^{1-\eps} \le n\le N$, we construct  RIP matrices
of order $k\ge n^{1/2+\eps}$ and constant $\delta=n^{-\eps}$.
This overcomes the natural barrier $k=O(n^{1/2})$ for proofs based on
small coherence, which are used in all previous explicit constructions of RIP
matrices.
Key ingredients in our proof are new
estimates for sumsets in product sets and for exponential sums
with the products of sets possessing special additive structure.
We also give a construction
of sets of $n$ complex numbers whose $k$-th moments are uniformly small
for $1\le k\le N$ (Tur\'an's power sum problem), which improves upon 
known explicit constructions when
$(\log N)^{1+o(1)} \le n\le (\log N)^{4+o(1)}$.  
This latter construction produces elementary explicit examples of $n \times N$
matrices that satisfy RIP and
whose columns constitute a new spherical code; for those problems the 
parameters closely match those of existing constructions
in the range $(\log N)^{1+o(1)}\le n\le (\log N)^{5/2+o(1)}$. 
\end{abstract}

\maketitle

%
\section{Introduction}
%

Suppose $1\le k \le  n\le N$ and $0<\delta<1$. A `signal' $\xx =
(x_j)_{j=1}^N \in \CC^N$ is said to be $k$-sparse if $\xx$
has at most $k$ nonzero coordinates. An $n \times N$   matrix
$\Phi$ is said to satisfy the
 Restricted Isometry Property (RIP) of order $k$ with constant
 $\delta$ if, for all $k$-sparse vectors $\xx$, we have
\begin{equation}\label{eq: RIP}
(1 - \delta)\|\xx\|_2^2 \le \|\Phi \xx\|_2^2 \le (1 + \delta)
\|\xx\|_2^2.
\end{equation}
While most authors work with real signals and matrices, in this paper
we work with complex matrices for convenience.  Given a complex 
matrix $\Phi$ satisfying \eqref{eq: RIP},  the $2n \times 2N$ real 
matrix $\Phi'$, formed by replacing each
element $a+ib$ of $\Phi$ by the $2 \times 2$ matrix 
$(\begin{smallmatrix} a&b\\ -b&a \end{smallmatrix})$, also
satisfies \eqref{eq: RIP} with the same parameters $k,\delta$.

We know from Cand\`es, Romberg and Tao that
matrices satisfying RIP have application to sparse signal recovery
(see \cite{Ca, CRT,CT}).  A variant of RIP (with the $\ell_2$ norm in
\eqref{eq: RIP} replaced by the $\ell_1$ norm) is also useful for such
problems \cite{BGIKS}.  A weak form of RIP, where \eqref{eq: RIP} holds
for most $k$-sparse $\xx$ (called Statistical RIP) is studied in \cite{GH}.
 Other applications of RIP matrices may be found in 
\cite{LT, NR}.

Given $n,N,\delta$, we wish to find $n \times N$ RIP matrices of
order $k$ with constant $\delta$, and with $k$ as large as
possible.  If the entries of $\Phi$ are independent Bernoulli
random variables with values $\pm 1/\sqrt{n}$, then with high
probability, $\Phi$ will have the required properties 
for\footnote{For convenience, we utilize the Vinogradov notation
$a\ll b$, which means $a=O(b)$, and the Hardy notation $a \order b$,
which means $b\ll a\ll b$.}
\be\label{kP} 
k \order \delta\frac{n}{\log (2N/n)}. 
\ee 
See \cite{CRT,MPT}; also \cite{BDDW} for a proof based on the
Johnson-Lindenstrauss lemma \cite{JL}.  The first result of
similar type for these matrices is due to Kashin \cite{Ka2}.
See also \cite{CT2,RV} for RIP matrices  with rows randomly  selected 
from the rows of a discrete Fourier transform matrix and for other 
random constructions of RIP matrices.  The parameter $k$ cannot be taken larger; in fact
\[
 k \ll \delta\frac{n}{\log (2N/n)}
\]
for every RIP matrix \cite{NT}.

It is an open problem to find good \emph{explicit}
constructions of RIP matrices; see T. Tao's Weblog \cite{Tao} for a
discussion of the problem.   We mention here that
all known explicit examples of RIP matrices are
based on constructions of systems of unit vectors (the columns of
the matrix)  with  small \emph{coherence}.


The \textit{coherence parameter} $\mu$ of a  collection of unit
vectors $\{\uu_1,\dots, \uu_N\}\subset \CC^n$ is defined by
\be\label{innerprod} 
\mu := \max_{r\ne s} |\langle \uu_r, \uu_s \rangle|. 
\ee 
Matrices whose columns are unit vectors with small
coherence are connected to a number of well-known problems, a few
of which we describe below.  Systems of vectors with small coherence are
also known as \emph{spherical codes}.
 Some other applications of
matrices with small coherence may be found in \cite{DET, GMS, Liv}.

\begin{prop}\label{cohRIP}
Suppose that $\uu_1,\dots,\uu_N$ are the columns of a matrix $\Phi$
and have coherence $\mu$.  Then $\Phi$ satisfies RIP of order $k$ with constant $\delta= (k-1)\mu$.
\end{prop}
\begin{proof}  For any $k$-sparse vector $\xx$,
\begin{align*}
|\|\Phi \xx\|_2^2-\|\xx\|_2^2| &\le 2\sum_{r<s}
|x_rx_s\langle \uu_r,\uu_s \rangle|\\
&\le \mu ((\sum |x_j|)^2 - \|\xx\|_2^2) \le (k-1)\mu \|\xx\|_2^2.\qedhere
\end{align*} 
\end{proof}

All explicit constructions of matrices with small coherence are based on
number theory.  There are many constructions producing matrices with
\be\label{expl_prev}
\mu \ll \frac{\log N}{\sqrt n\log n}.
\ee
In particular, such examples have been constructed by
Kashin \cite{Ka}, Alon, Goldreich, H{\aa}stad and Peralta \cite{Alon},
 DeVore \cite{D}, and Nelson and Temlyakov \cite{NT}.
By Proposition \ref{cohRIP},
these matrices satisfy RIP with constant $\delta$ and order
\be\label{kD} 
k \order \delta\frac{\sqrt{n}\log n}{\log N}. 
\ee
It follows from random constructions of Erd\H os and R\'enyi for Tur\'an's
problem (see Proposition \ref{Turancoh} and \eqref{upE-R} below) that for
any $n,N$ there are vectors with coherence
\[
 \mu \ll \sqrt{\frac{\log N}{n}}.
\]
By contrast, there is a universal lower bound
\be\label{mulower}
\mu \gg \Big( \frac{\log N}{n\log(n/\log N)} \Big)^{1/2}
\ge \frac{1}{\sqrt{n}},
\ee
valid for $2\log N\le n \le N/2$ and all $\Phi$, due to
Levenshtein \cite{Lev} (see also \cite{Glu} and \cite{NT}). 
Therefore, 
by estimating RIP parameters in terms of
the coherence parameter we cannot construct $n \times N$ RIP
matrices of order larger than $\sqrt{n}$ and constant $\delta<1$.

Using methods of additive
combinatorics, we construct RIP matrices of order $k$ with
$n=o(k^2)$.

\begin{thm}\label{thmRIP}
There is an effective constant $\eps_0>0$ and an explicit number
$n_0$ such that for any positive integers $n\geq n_0$ and 
$n \le N\le  n^{1+\eps_0}$, there is an explicit $n \times N$ RIP matrix of order
$\lfloor n^{\frac12+\eps_0}\rfloor$ with constant $n^{-\eps_0}$.
\end{thm}

\begin{remark}
For application to sparse signal recovery, it is sufficient to take
fixed $\delta < \sqrt{2}-1$ \cite{Ca}, and one needs an upper bound on $n$
in terms of $k,N$.  By Theorem \ref{thmRIP}, for some $\eps_0'>0$, large $N$
and $N^{1/2-\eps_0'} \le k \le N^{1/2+\eps_0'}$, we construct explicit RIP
matrices with $n\le k^{2-\eps_0'}$.
\end{remark}

The proof of Theorem \ref{thmRIP} uses a result on additive energy
of sets (Corollary \ref{measener}, Theorem \ref{measener2}),
estimates for sizes of sumsets in product sets (Theorem
\ref{sumset}), and bounds for exponential sums over products of
sets possessing special additive structure (Lemma \ref{2ndcase}).

\bigskip

We now return to the problem of constructing matrices with small coherence.
By \eqref{mulower}, the bound \eqref{expl_prev} cannot be improved if
$\log n \gg \log N$, but there is a gap between bounds \eqref{mulower}
and \eqref{expl_prev} when $\log n = o(\log N)$.  For example,
\eqref{expl_prev} is nontrivial only for $n\gg (\log N/\log\log N)^2$.
Of particular interest in coding theory is the range $n=O(\log^C N)$ for
fixed $C$, where there have been some improvements made to \eqref{expl_prev}.
A construction obtained by concatenating algebraic-geometric codes
with Hadamard codes (see e.g. \cite[Corollary 3]{GS} and Section 3 
of \cite{BenTa}) produces  matrices with coherence
\be\label{AGH}
\mu \ll \pfrac{\log N}{n\log(n/\log N)}^{1/3},
\ee
which is nontrivial for $n \gg \log N$, and is better than
\eqref{expl_prev} when $\log N \ll n \ll (\frac{\log N}{\log\log N})^4$.
In the range $(\frac{\log N}{\log\log N})^{5/2} \ll n \ll (\frac{\log N}{\log\log N})^5$, Ben-Aroya and Ta-Shma \cite{BenTa} improved both
\eqref{expl_prev} and \eqref{AGH} by constructing binary codes
(vectors with entries $\pm 1/\sqrt{n}$) with coherence
\be\label{BenTa}
 \mu \ll \pfrac{\log N}{n^{4/5} \log\log N}^{1/2}.
\ee

In this paper, we introduce very elementary constructions of matrices with
coherence which matches (up to a $\log\log N$ factor) the bound \eqref{AGH}.
Our constructions, which are based on a method of
Ajtai, Iwaniec, Koml\'os, Pintz and Szemer\'edi \cite{5}, 
have the added utility of applying to Tur\'an's power-sum
problem and to the problem of finding thin sets with small Fourier coefficients.
For the last two problems, our construction gives better estimates
than existing explicit constructions in certain ranges of the parameters.

Roughly speaking, a set with small Fourier coefficients can be used to 
construct a set of numbers for Tur\'an's problem, and a set of numbers in Tur\'an's
problem can be used to produce a matrix with small coherence.  
This is made precise below.

We next describe the problem of explicitly constructing thin sets with small
Fourier coefficients.  If $N$ is a positive integer and $S$ is a set (or multiset) of residues modulo $N$, we let
$$
f_S(k) = \sum_{s\in S} e^{2\pi i ks/N}
$$
and
$$
|f_S| := \frac{1}{|S|} \max_{1\le k\le N-1} |f_S(k)|.
$$
Given $N$, we wish to find a small set $S$ with $| f_S |$ also small.

Tur\'an's problem \cite{Tu} concerns the estimation of the function
$$
T(n,N) = \min_{|z_1|=\dots=|z_n|=1} M_N(\zz), \quad 
M_N(\zz) := \max_{k=1,\dots,N} \left|\sum_{j=1}^n z_j^k\right|.
$$
where $n, N$ are positive integers. 
There is a vast literature related to Tur\'an's problem;
see, e.g., \cite{An1}, \cite{An2}, \cite{Mon} (chapter 5),
\cite{St1}, \cite{St2}.  

If $S=\{t_1,\dots,t_n\}$ is a multiset of integers modulo $N$ and $z_j=e^{2\pi i t_j/N}$ for $1\le j\le n$, we see that
\be\label{TfS}
T(n,N-1) \le M_{N-1}(\zz) \le n|f_S|.
\ee
We also have the following easy connection between Tur\'an's problem and
coherence.
\begin{prop}\label{Turancoh}
Given any vector $\zz=(z_1,\ldots,z_n)$ with $|z_j|=1$ for all $j$,
the coherence $\mu$ of the $n\times N$ matrix with the columns
\be\label{ukzn}
\uu_k^{-1/2}(z_1^{k-1},\dots,z_n^{k-1})^T,\quad k=1,\dots,N
\ee
satisfies $\mu = n^{-1} M_{N-1}(\zz)$.
\end{prop}

Combining \eqref{TfS} and Proposition \ref{Turancoh}, for any multiset $S$
of residues modulo $N$, the vectors \eqref{ukzn} satisfy
\be\label{mufS}
\mu \le | f_S |.
\ee

A corollary of a character sum estimate of  Katz \cite{Kat} (see also \cite{PS})
shows\footnote{Here we take $N=p^d-2$, where $p$ is prime, $p\approx 
((d-1)/\mu)^2$ and $((d-1)\mu^{-1})^{2d} \approx N$.  Let 
$F=\mathbb{F}_{p^d}$.  The group of characters on $F$ is a cyclic group of 
order $N+1$ with generator $\chi_1$.   For any $x\in F\setminus\{0\}$ write
$\chi_1(x) = e(t_x/N)$. Let $x$ be an element of $F$ not contained in any 
proper subfield of $F$ and take $T=\{t_{x+j}: j=0,\dots,p-1\}$.  Then $|T|=p$,
and $|f_T|\le(d-1)\sqrt p$ by \cite{Kat}.
}
that for certain $N$ and $1/N \le \mu \le 1$, there are  
(explicitly defined) sets $T$ of residues modulo $N$ so that
\be\label{Katz}
|f_T| \le \mu, \qquad |T| = O\pfrac{\log^2 N}{\mu^2(\log\log N+\log(1/\mu))}.
\ee

An application of Dirichlet's approximation theorem
shows that a set $S$ with $|S| < \log N$ must have $|f_S| \gg 1$.
In \cite{5}, sets which are not much larger are explicitly constructed
so that $|f_S|$ is small.  Specifically, by \cite[(1),(2)]{5}, 
for each prime\footnote{A corresponding result when $N$ is composite is given 
in \cite{RSW}.}
$N$ there is a set $S$ with $|S|=O(\log N (\log^* N)^{13\log^* N})$ and
$$
|f_S| = O(1/\log^* N),
$$
where $\log^* N$ is the integer $k$ so that the $k$-th iterate of the logarithm
of $N$ lies in $[1,e)$.  The proof uses an iterative procedure.  By modifying
this procedure, and truncating after two steps, we prove the following.
To state our results, for brevity write
$$
L_1=\log N, \qquad L_2 = \log\log N, \qquad L_3 = \log\log\log N.
$$

\begin{thm}\label{thm1}
For sufficiently large prime $N$ and $\mu$ such that
\be\label{KN}
\frac{L_2^4}{L_1} \le \mu < 1, \qquad 1/\mu \in \mathbb{N},
\ee
a set $S$ of residues modulo $N$ can be explicitly constructed so that
$$
|f_S| \le \mu, \quad \text{and} \quad 
|S| =O\pfrac{L_1 L_2\log(2/\mu)}{\mu^4(L_3+\log(1/\mu))} = 
O\pfrac{L_1L_2}{\mu^4}.
$$
\end{thm}

\begin{remark}
The method from \cite{5}, if applied without modification (with two iterations
of the basic lemma), produces a conclusion in Theorem \ref{thm1} with
$$
|S| =O\pfrac{L_1 L_2}{\mu^8 L_3}.
$$
\end{remark}

\begin{remark}
 The bound on $|S|$ in Theorem \ref{thm1} is better than \eqref{Katz}
for $\mu \gg L_1^{-1/2}L_2$.
\end{remark}

Together, the construction for Theorem \ref{thm1} and \eqref{TfS} give
explicit sets $\zz$ for Tur\'an's problem.  By further modifying the construction,
we can do better.

\begin{thm}\label{thm2}
For sufficiently large positive integer $N$ and $\mu$ such that
\be\label{KN2}
\frac{L_2^3}{L_1} \le \mu < 1,
\ee
a multiset $\zz=\{z_1,\dots,z_n\}$ such that $|z_1|=\dots=|z_n|=1$,
can be explicitly constructed so that
$$
M_N(\zz) \le \mu n, \qquad
n =O\pfrac{L_1 L_2\log(2/\mu)}{\mu^3(L_3+\log(1/\mu))}=O\pfrac{L_1L_2}{\mu^3}.
$$
\end{thm}

To put Theorem \ref{thm2} in context, we briefly review what is known
about $T(n,N)$.  
P. Erd\H os and A. R\'enyi \cite{ER} used
probabilistic methods to prove an upper estimate
\be\label{upE-R}
T(n,N) \le (6n\log(N+1))^{1/2}.
\ee
Using the character sum bound of Katz \cite{Kat}, J. Andersson \cite{An3}
gave explicit examples of sets $\zz$ which give
\be\label{upAnd}
T(n,N) \le M_N(\zz) \ll \frac{\sqrt n\log N}{\log n}
\ee
 One can see that
\eqref{upAnd} supersedes \eqref{upE-R} for $\log N\ll\log^2 n$.  Also,
combining \eqref{upAnd} with Proposition \ref{Turancoh} provides 
yet another construction
of matrices with coherence satisfying \eqref{expl_prev}.
On the other hand, by \eqref{mulower} and Proposition \ref{Turancoh}, 
we have the lower estimate
\[
 T(n,N) \gg \pfrac{n \log N}{\log (n/\log N)}^{1/2} \gg n^{1/2}
\qquad (2\log N \le n\le N/2).
\]

By comparison, the constructions in Theorem \ref{thm2} are better than
\eqref{upAnd} in the range $n\ll L_1^4/L_2^8$, that is, throughout the range
\eqref{KN2} (our constructions require $n$ to be prime, however).

The constructions in Theorem \ref{thm2} also produce, by 
Proposition \ref{Turancoh},
 explicit examples of matrices with coherence
\[
 \mu  \ll \pfrac{L_1L_2}{n}^{1/3},
\]
which is close to the bound \eqref{AGH}.  By Proposition \ref{cohRIP}, these
matrices satisfy RIP with constant $\delta$ and order
\[
 k \gg \delta \pfrac{n}{L_1L_2}^{1/3}.
\]


We prove Theorem \ref{thmRIP} in Sections \ref{Sec2}--\ref{Sec7},
Theorem \ref{thm1} in Section \ref{Sec8} and Theorem \ref{thm2} in Section
\ref{Sec9}.

%
%
\section{Construction of the matrix in Theorem~\ref{thmRIP}}\label{Sec2}
%
%

We fix a large even number $m$. A value of $m$ can be specified;
it depends on the constant $c_0$ in an estimate from additive combinatorics
(Proposition \ref{TheoremC}, Section \ref{Sec4}).
Also, the value $m$ can be reduced if one proves a better version of the
 Balog--Szemer\'edi--Gowers lemma (Lemma \ref{BSG} below).

For sufficiently large $n$ we take the largest prime $p\leq n$, which
satisfies $p\ge n/2$ by Bertrand's postulate.
By $\Fp$ we denote the field of the residues modulo $p$,
and let $\Fps=\Fp\setminus\{0\}$.
For $x\in\Fp$, let $e_p(x)=e^{2\pi i x/p}$.
We construct an appropriate $p \times N$
matrix $\Phi_p$ with columns $\uu_{a,b}, a\in\CA\subset\Fp, b\in\CB\subset\Fp$
where
$$
\uu_{a,b}=\frac1{\sqrt p}(e_p(ax^2+bx))_{x\in\Fp}
$$
and the sets $\CA, \CB$ will be defined below. Notice that the matrix
$\Phi_p$ can be extended to a $n \times N$ matrix $\Phi$ by adding
$n-p$ zero rows. Clearly, the matrices $\Phi_p$ and $\Phi$ have
the same RIP parameters.

We take
\be\label{Adef}
\a=\frac{1}{8m^2}, \quad L=\lfloor p^{\a}\rfloor,\quad U=L^{4m-1},
\qquad \CA=\{x^2+Ux: 1\le x\le L\}.
\ee
To define the set $\CB$, we take
\[
\beta=\alpha/2=1/(16m^2),\quad r=\left\lfloor\frac{\beta\log p}
{\log 2}\right\rfloor, \quad M=2^{(1/\beta)-1}=2^{16m^2-1},
\]
and let
$$
\CB=\left\{\sum_{j=1}^r x_j(2M)^{j-1}:\,x_1,\dots,x_r\in\{0,\dots,M-1\}
\right\}.
$$
We notice that all elements of $\CB$ are at most $p/2$, and
\be\label{Bbound}
|\CB|\asymp p^{1-\beta}.
\ee
It follows from \eqref{Adef} and \eqref{Bbound} that
$$
|\CA||\CB|\asymp p^{1+\beta}\asymp n^{1+\beta}.
$$
For $n \le N \le n^{1+\beta/2}$, take $\Phi$ to be the
matrix formed by the first $N$ columns of $\Phi_p$, padded with
$n-p$ rows of zeros.

In the next four sections, we show that $\Phi$ has the required properties
for Theorem \ref{thmRIP}.  First, in Section \ref{section:flatrip}, we show
that  in \eqref{eq: RIP} we need only consider vectors $\xx$ whose components 
are 0 or 1 (emph{flat} vectors).  We prove the following.

\begin{lem}\label{flat2} Let $k\ge2^{10}$ and $s$ be a positive integer.
Assume that the coherence parameter of the matrix $\Phi$ is $\mu\le1/k$.
Also, assume that for some $\delta\ge0$ and any disjoint
$J_1,J_2\subset\{1,\dots,N\}$ with $|J_1|\le k, |J_2|\le k$ we have
$$\left|\left\langle\sum_{j\in J_1}\uu_j,\sum_{j\in J_2}\uu_j\right
\rangle\right| \le\delta k.
$$
Then $\Phi$ satisfies the RIP of order $2sk$ with constant
$44s\sqrt{\delta}\log k$.
\end{lem}

Our main lemma concerns showing RIP with flat vectors and 
order $k=\lfloor \sqrt{p} \rfloor$.  We prove the required estimates for
matrices formed from more general sets $\CA$ and $\CB$ having
certain additive properties.  Namely, 
let $m\in2\NN$ and $0 < \alpha < 0.01$.
Assume that
\be\label{Aupper}
|\CA|\le p^\alpha
\ee
and, for $a\in \CA$ and $a_1,\dots,a_{2m}\in\CA\setminus\{a\}$,
\be\label{Asparse}
\sum_{j=1}^{m}\frac{1}{a-a_j}=\sum_{j=m+1}^{2m} \frac{1}{a-a_j}\;
\implies \; (a_1,\ldots,a_m)\text{ is a permutation of }(a_{m+1},\ldots,a_{2m}).
\ee
Here we write $1/x$ for the multiplicative inverse of $x\in\Fp$.
We will consider the sets $\CB$ satisfying
\be\label{Benergy}
\forall S\subset\CB\quad\text{if}\quad|S|\ge p^{1/3}\quad\text{then}
\quad E(S,S)\le p^{-\gamma}|S|^3
\ee
with some $\gamma>0$, where $E(S,S)$ is the number of solutions of
$s_1+s_2=s_3+s_4$ with each $s_i\in S$.

\begin{lem}\label{main} Let $m\in2\NN$, $\alpha\in(0,0.01)$, $0<\gamma \le
\min(\a,\frac{1}{3m})$, $p$ sufficiently large in terms of  $m,\alpha,\gamma$, $\CA$ satisfies \eqref{Aupper} and \eqref{Asparse},
and $\CB$ satisfies $\eqref{Benergy}$.  Then for any disjoint sets
$\Omega_1,\Omega_2\subset\CA\times\CB$ such that $|\Omega_1|\le\sqrt p$,
$|\Omega_2|\le\sqrt p$, the inequality
$$
\left|\sum_{(a_1,b_1)\in \Omega_1}\sum_{(a_2,b_2)\in \Omega_2}
\left\langle \uu_{a_1,b_1}, \uu_{a_2,b_2}\right\rangle\right|
\leq p^{1/2-\eps_1}$$
holds where $\eps_1=\cml \g/20-43 \a/m$.
\end{lem}

The proof of Lemma \ref{main} is quite involved, and will be
handled in three subsequent sections.  We next demonstate
how Theorem \ref{thmRIP} may be deduced from it.

We first prove \eqref{Asparse} for the specific set $\CA$ defined in
\eqref{Adef}, provided that $p>(2m)^{8m^2}$
(and thus $L \ge 2m$).  We have to show that for any distinct
$x,x_1\dots,x_n\in\{1,\dots,L\}$ and any nonzero integers
$\lambda_1,\dots,\lambda_n$ such that $n\ge 2m$ and
$|\lambda_1|+\cdots+|\lambda_n|\le 2m,$
the sum
$$
V=\sum_{j=1}^n\frac{\lambda_j}{(x-x_j)(x+x_j+U)}
$$
is a nonzero element of $\Fp$. However, we will treat $V$
as a rational number. Denote
$$D_1=\prod_{j=1}^n(x-x_j),\quad D_2=\prod_{j=1}^n(x+x_j+U).$$
So,
\be\label{sumfract}
D_1D_2V=\sum_{j=1}^n\frac{\lambda_jD_1}{x-x_j}\frac{D_2}{x+x_j+U}.
\ee
All summands in the right-hand side of \eqref{sumfract} but the first one
are divisible by $x+x_1+U$. For the first summand we have
$$
\frac{\lambda_1D_1}{x-x_1}\frac{D_2}{x+x_1+U}\equiv V_1 \; (\bmod\; x_0+x_1+U),
$$
where
$$
V_1=\lambda_1\prod_{j=2}^n(x-x_j)\prod_{j=2}^n(x_j-x_1).
$$
We have
$$
|V_1|\le 2mL^{2n-2}\le 2mL^{4m-2} \le L^{4m-1} = U < U+x_0+x_1.
$$
This shows that $V_1\ne 0 \; (\bmod\; x_0+x_1+U)$. Therefore, $V\neq0$.
By assumption, $p\nmid D_1$, and
\[
|D_2V|\le 2m(U+2L)^n/U \le 4m U^{2m-1} \le U^{2m} <p. \qedhere
\]
Hence $p\nmid D_1 D_2 V$, as desired.

Condition~\eqref{Benergy} is satisfied due to Corollary~\ref{sumset2}
of Section \ref{Sec5} with $\gamma=\beta/50$.
 If $m>86000\cml^{-1}$ then Lemma~\ref{main}
gives a nontrivial estimate with $\eps_1>0$.
Thus, $\Phi_p$ satisfies the conditions of Corollary
\ref{flat2} with
$k=\lfloor \sqrt{p}\rfloor \ge \sqrt{n/2}$ and
$\delta=p^{-\eps_1} \le (n/2)^{-\eps_1}$ (using $p\ge 0.9n$
for large $n$, which follows from the prime number theorem).
Let $\eps_0=\eps_1/5$.  Let $n\le N\le n^{1+\eps_0}$,
and let $\Phi$ be the $n \times N$ matrix formed by taking the
first $N$ columns of $\Phi_p$, then adding $n-p$ rows of zeros.
Clearly, $\Phi$ satisfies the conditions of Corollary \ref{flat2} with
the same  parameters as $\Phi_p$.
By Lemma \ref{flat2} with $s=\lfloor p^{\eps_1/4} \rfloor$,
Theorem~\ref{thmRIP} follows.

In Section \ref{Sec4} we introduce some notation and recall standard
estimates in additive combinatorics, which will be applied to subsets
of $\CB$.  Section \ref{Sec5} is devoted to the sumset theory of $\CB$,
from which we deduce \eqref{Benergy}.  
The completion of the proof of
Lemma \ref{main} is in Section \ref{Sec7}.  We give some preliminaries
here.

It is easy to see that for a fixed $a$ the vectors
$\{u_{a,b}:\,b\in\Fp\}$ form an orthogonal system.
Using a well-known formula for Gauss sums $\sum_{x\in\Fp}e_p(dx^2)$
(see, for example, \cite{IR}, Proposition~6.31), we have
for $a_1\neq a_2$ the equality
\begin{align*}
\langle \uu_{a_1,b_1},\uu_{a_2,b_2}\rangle &=p^{-1}
e_p\left(-\frac{(b_1-b_2)^2}{4(a_1-a_2)}\right)\sum_{x\in\Fp}
e_p((a_1-a_2)x^2) \\ &=\frac{\sigma_p}{\sqrt p}
\left(\frac{a_1-a_2}p\right)e_p\left(-\frac{(b_1-b_2)^2}
{4(a_1-a_2)}\right),
\end{align*}
where $(\frac dp)$ is the Legendre 
symbol\footnote{for $d\in\Fps$, we have $(\frac dp)=1$
if the congruence $x^2\equiv d\pmod{p}$ has a solution, and $(\frac dp)=-1$
otherwise.}, and $\sigma_p=1$ or $i$
according as $p\equiv 1$ or $3\pmod 4$.  We remark that
there is no analogous formula for exponential sums $\sum_{x\in\Fp}
e_p(F(x))$ when $F$ is a polynomial of degree $\ge 3$.
Consequently, the assertion of Lemma~\ref{main} can be
rewritten as
\be\label{mainineq}
\left|\sum_{(a_1,b_1)\in \Omega_1}\sum_{(a_2,b_2)\in \Omega_2}
\left(\frac{a_1-a_2}p\right) e_p\left(\frac{(b_1-b_2)^2}
{4(a_1-a_2)}\right)\right|\leq p^{1-\eps_1},
\ee
where the 
summands with $a_1=a_2$ are excluded from the summation. 
We next break $\Omega_1, \Omega_2$ into \emph{balanced} sets.
For $a\in\CA$ and $i=1,2$, let
$$
\Omega_i(a)=\{b\in\CB:\,(a,b)\in\Omega_i\}.
$$
To prove \eqref{mainineq} it is enough to show that
\be\label{balanced}
|S(A_1,A_2)| \le p^{1-1.1\eps_1}, \quad S(A_1,A_2)=
\sum_{\substack{a_1\in A_1,\\a_2\in A_2}}
\sum_{\substack{b_1\in\Omega_1(a_1),\\b_2\in\Omega_2(a_2)}}
\left(\frac{a_1-a_2}p\right) e_p\left(\frac{(b_1-b_2)^2}
{4(a_1-a_2)}\right),
\ee
whenever $M_1, M_2$ are powers of two and, for $i=1,2$ and for any $a_i\in A_i$,
\be\label{balanceM}
M_i/2\le|\Omega_i(a_i)|<M_i, \qquad |A_i| M_i \le 2\sqrt{p}.
\ee
Indeed, there are $O(\log^2 p)$ choices for $M_1,M_2$.
To prove the cancellation in \eqref{balanced}, we basically split into
two cases: (i) some $B'=\Omega_i(a_j)$ has additive structure (that is,
$E(B',B')$ is large), where the cancellation comes from the sum over
$b_1,b_2$ (with $a_1,a_2$ fixed), and (ii) when
$B'$ does not have additive structure, in
which case one gets dispersion of the phases from the dilation weights
 $1/(a_1-a_2)$ (taking a large moment and using \eqref{Asparse}).
Incidentally, oscillations of the factor $(\frac{a_1-a_2}{p})$ 
play no role in the argument.

%
%
\section{The Flat-RIP property}\label{section:flatrip}
%
%

Let $\uu_1,\ldots,\uu_N$ be the columns of an $n \times N$ matrix $\Phi$.
Suppose that for every $j$, $\|\uu_j\|_2=1$. We say that $\Phi$
satisfies the flat RIP of order $k$ with constant $\delta$
if for any disjoint $J_1,J_2\subset\{1,\dots,N\}$
with $|J_1|\le k, |J_2|\le k$ we have
\begin{equation} \label{eq: flatRIP}
\left|\left\langle\sum_{j\in J_1}\uu_j,\sum_{j\in J_2}\uu_j\right\rangle\right|
\le\delta(|J_1||J_2|)^{1/2}.
\end{equation}

For technical reasons, it is more convenient to work with the flat-RIP than
with the RIP. However, flat-RIP implies RIP with an
increase in $\delta$.  The flat-RIP property is closely related to the
property that \eqref{eq: RIP} holds for any $\xx$ with entries which are
zero or one and at most $k$ ones (see the calculation at the end of this
section).

\begin{lem}\label{flat1} Let $k\ge2^{10}$ and $s$ be a positive integer.
Suppose that $\Phi$ satisfies flat-RIP of order $k$ with constant $\delta$.
Then $\Phi$ satisfies RIP of order $2sk$ with constant
$44s\delta\log k$.
\end{lem}

\begin{proof}  First, by a convexity-type argument and our assumption,
\begin{equation} \label{eq: [0,1]}
\left|\left\langle\sum_{j\in J_1}x_j\uu_j,\sum_{j\in J_2}y_j\uu_j\right
\rangle\right|\le\delta(|J_1||J_2|)^{1/2}
\end{equation}
provided that $|J_1|\le k, |J_2|\le k$, $0\le x_j,y_j\le1$ for all $j$.
Next, suppose $|J_1|\le k, |J_2|\le k$, and $0\le x_j,y_j$ for all $j$.
Without loss of generality assume that
$\|\xx\|_2=\|\yy\|_2=1$, where $\| \cdot\|_2$ denotes the $l_2$ norm.
For a positive integer $\nu$ let
$$
J_{1,\nu}=\{j\in J_1:\,2^{-\nu}<x_j\le2^{1-\nu}\},\quad
J_{2,\nu}=\{j\in J_2:\,2^{-\nu}<y_j\le2^{1-\nu}\}.
$$
Observe that
\begin{equation} \label{estJ}
\sum_\nu 4^{-\nu}|J_{1,\nu}|\le1,\quad\sum_\nu 4^{-\nu}|J_{2,\nu}|\le1.
\end{equation}
Applying \eqref{eq: [0,1]} to sets $J_{1,\nu}, J_{2,\nu}$, we get
\begin{align*}
\left|\left\langle\sum_{j\in J_1}x_j\uu_j,\sum_{j\in J_2}y_j\uu_j\right
\rangle\right|
&\le\sum_{\nu_1,\nu_2}\left|\left\langle\sum_{j\in J_{1,\nu_1}
}x_j\uu_j,
\sum_{j\in J_{2,\nu_2}}y_j\uu_j\right\rangle\right|\\
&\le\sum_{\nu_1,\nu_2}2^{2-\nu_1-\nu_2}\delta(|J_{1,\nu_1}||J_{2,\nu_2}|)^{1/2}
\\
&=4\delta\sum_{\nu}2^{-\nu}|J_{1,\nu}|^{1/2}
\sum_{\nu}2^{-\nu}|J_{2,\nu}|^{1/2}.
\end{align*}
Let $t=\lfloor 3+\log k/(2\log 2)\rfloor$. By the \CS inequality we infer
that
\begin{align*}
\sum_{\nu}2^{-\nu}|J_{1,\nu}|^{1/2}&\le
\sum_{\nu=1}^t 2^{-\nu}|J_{1,\nu}|^{1/2}
+ \sum_{\nu=t+1}^\infty 2^{-\nu}|J_{1,\nu}|^{1/2}\\
&\le t^{1/2}\left(\sum_{\nu=1}^t 4^{-\nu}|J_{1,\nu}|\right)^{1/2}
+ \sum_{\nu=t+1}^\infty 2^{-\nu}k^{1/2}\le t^{1/2}+\frac14.
\end{align*}
Similarly,
$$\sum_{\nu}2^{-\nu}|J_{2,\nu}|^{1/2}\le t^{1/2}+\frac14.$$
Therefore,
\be\label{l2pos}
\left|\left\langle\sum_{j\in J_1}x_j\uu_j,\sum_{j\in J_2}y_j\uu_j\right
\rangle\right|\le4\delta\left(t^{1/2}+\frac14\right)^2
\le 5.5 \delta \log k.
\ee

For the next step, suppose $x_j,y_j$ take arbitrary
complex values, $|J_1|\le sk$ and $|J_2|\le sk$.
We partition $J_1$ and $J_2$ into $s$ subsets of cardinality at most $k$
each:
$J_1=\cup_{\mu=1}^s J_{1,\mu}$,  $J_2=\cup_{\mu=1}^s J_{2,\mu}.$
Next, for any $j$ we have
$$
x_j=\sum_{\nu=1}^4 x_{j,\nu}i^{\nu},\quad
y_j=\sum_{\nu=1}^4 y_{j,\nu}i^{\nu}, \quad
|x_j|^2=\sum_{\nu=1}^4 x_{j,\nu}^2,\quad
|y_j|^2=\sum_{\nu=1}^4 y_{j,\nu}^2,
$$
where $x_{j,\nu}, y_{j,\nu}$ are non-negative.
By \eqref{l2pos} and the \CS inequality,
\be\label{l2com}\begin{split}
\left|\left\langle\sum_{j\in J_1}x_j\uu_j,\sum_{j\in J_2}y_j\uu_j\right
\rangle\right|
&\le\sum_{\mu_1=1}^s\sum_{\nu_1=1}^4\sum_{\mu_2=1}^s\sum_{\nu_2=1}^4
\left|\left\langle\sum_{j\in J_{1,\mu_1}}x_{j,\nu_1}\uu_j,
\sum_{j\in J_{2,\mu_2}}y_{j,\nu_2}\uu_j\right\rangle\right| \\
&\le\sum_{\mu_1, \nu_1, \mu_2, \nu_2}
5.5\delta(\log k)\left(\sum_{j\in J_{1,\mu_1}}x_{j,\nu_1}^2\right)^{1/2}
\left(\sum_{j\in J_{2,\mu_2}}y_{j,\nu_2}^2\right)^{1/2}\\
&\le22s\delta\|\xx\|_2\|\yy\|_2\log k.
\end{split}\ee

To complete the proof of the lemma
assume $N\ge 2sk$ and consider a vector
$\xx=\sum_{j\in J}x_je_j$ with $\| \xx \|_2=1$ and $|J|=2sk,$
where $(e_1,\ldots,e_N)$ is the standard basis of $\mathbb{C}^N$.
Take arbitrary partitions of $J$ into two sets $J_1, J_2$ of cardinality $sk$
each. By \eqref{l2com}, we have
\begin{align*}
\left| \| \Phi \xx \|_2^2 - \| \xx \|_2^2 \right| &=
\left|\sum_{j_1,j_2\in J, j_1\ne j_2}
\langle x_{j_1}\uu_{j_1},x_{j_2}\uu_{j_2}\rangle\right|\\
&=\binom{2sk-2}{sk-1}^{-1}\left|\sum_{J_1,J_2}
\left\langle\sum_{j\in J_1}x_j\uu_j,
\sum_{j\in J_2}x_j\uu_j\right\rangle\right|\\
&\le\binom{2sk-2}{sk-1}^{-1}\sum_{J_1,J_2}
22s\delta(\log k)\left(\sum_{j\in J_1}|x_j|^2\right)^{1/2}
\left(\sum_{j\in J_2}|x_j|^2\right)^{1/2}\\
&\le\binom{2sk-2}{sk-1}^{-1}\sum_{J_1,J_2}11s\delta\|\xx\|_2^2\log k\\
&=\binom{2sk}{sk}\binom{2sk-2}{sk-1}^{-1}11s\delta\|\xx\|_2^2\log k
\le 44s\delta\|\xx\|_2^2\log k.
\qedhere
\end{align*}
\end{proof}

\begin{proof}[Proof of Lemma \ref{flat2}] For any disjoint $J_1,J_2\subset\{1,\dots,N\}$ with
$|J_1|\le k, |J_2|\le k$ we have
$$\left|\left\langle\sum_{j\in J_1}\uu_j,
\sum_{j\in J_2}\uu_j\right\rangle\right|
\le\min(\delta k,\mu|J_1||J_2|)\le\min(\delta k,|J_1||J_2|/k)
\le\sqrt{\delta|J_1||J_2|},$$
and it remains to apply Lemma~\ref{flat1}.
\end{proof}

\begin{remark} Using the assumptions of the Lemma \ref{flat2} directly
rather than reducing it to Lemma~\ref{flat1}, one can get
a better constant for RIP; 
However, we do not need a stronger version of the corollary for
 our purposes.
\end{remark}

%
%
\section{Some definitions and results from additive combinatorics}\label{Sec4}
%
%

For an (additive) abelian group $G$ we define the sum and the difference
of subsets $A,B\subset G$:
$$A+B = \{a+b:\,a\in A, b\in B\},\quad A-B = \{a-b:\,a\in A, b\in B\}.$$
We denote $-A=\{-x: x\in A \}$.  If $A \subseteq G=\Fp$ and
$b\in \Fp$, write $bA=\{ ba: a\in A \}$.

Consider $G=\Fp$ and let $\CB\subset G$ be the set defined in Section
\ref{Sec2}.  There is a natural bijection $\Phi$ between $\CB$ and the cube
$\cC_{M,r}=\{0,\dots,M-1\}^r$ defined by
$\Phi(\sum_{j=1}^r x_j(2M)^{j-1})=(x_1,\dots,x_r).$  Moreover, it is trivial
that $b_1+b_2=b_3+b_4$  if and only if
$\Phi(b_1)+\Phi(b_2)=\Phi(b_3)+\Phi(b_4)$.  In the language
of additive combinatorics, $\Phi$ is a Freiman
isomorphism between $\CB$ and $\cC_{M,r}$.
Thus, $|B_1+B_2| = |\Phi(B_1)+\Phi(B_2)|$ for any $B_1\subseteq \CB$, $B_2
\subseteq \CB$.
The problem of the size of sumsets in $\cC_{M,r}$ will be investigated in the next section.

We will use the following lemma which is a particular case of Pl\"unecke --
Ruzsa estimates (\cite{TaoVu}, Exercise 6.5.15).

\begin{lem}\label{PlunRuz}
For any nonempty set $A\subset G$ we have $|A+A|\le|A-A|^2/|A|$.
\end{lem}

If $A,B\subset G$, we define the (additive) energy $E(A,B)$ of the sets $A$ and
$B$ as the number of solutions of the equation
$$a_1+b_1=a_2+b_2,\quad a_1,a_2\in A,\,b_1,b_2\in B.$$
Next, let $F\subset A\times B$. The $F$-restricted sum of $A$ and $B$
is defined as
$$A+_F B = \{a+b:\,a\in A, b\in B, (a,b)\in F\}.$$
Trivially $E(A,A)\le|A|^3.$
If $E(A,A)$ is close to $|A|^3$ then $A$ must have a special
additive structure.

\begin{lem}\label{enersum} (\cite{TaoVu}, Lemma~2.30) If $E(A,A)\ge|A|^3/K$
then there exists $F\subset A\times A$ such that $|F|\ge|A|^2/(2K)$
and
$|A+_F A|\le2K|A|$.
\end{lem}

The following lemma \cite{BG} is a version of the Balog--Szemer\'edi--Gowers
lemma which plays a very important role in additive combinatorics.

\begin{lem}\label{BSG} If $F\subset A\times A$, $|F|\ge|A|^2/L$
and $|A+_F A|\le L|A|$. Then there exists a set $A'\subset A$
such that $|A'|\ge|A|/(10L)$ and $|A'-A'|\le 10^4 L^9|A|$.
\end{lem}

Combining  Lemma~\ref{enersum} and Lemma~\ref{BSG} gives the following.
\begin{cor}\label{BSG2} If $E(A,A)\ge|A|^3/K$
then there exists a set $A'\subset A$ such that
$|A'|\ge|A|/(20K)$ and $|A'-A'|\le10^7K^9|A|$.
\end{cor}

For a function $f:\Fp\to\CC$ and a number $r\ge1$ we define the $L_r$ norm
of $f$:
$$
\|f\|_r=\Bigg(\sum_{x\in\Fp}|f(x)|^r\Bigg)^{1/r}.
$$
The additive convolution of two functions $f,g:\Fp\to\CC$
is defined as
$$f*g(x)=\sum_{y\in\Fp}f(y)g(x-y).$$

By $1_A$ we denote the indicator function
of the set $A$. With this notation,  we have
\be\label{EAB1A}
E(A,B)=E(A,-B)=\|1_A*1_B\|_2^2.
\ee

We say that a function $f:\Fp\to\RR_+$ is a probability measure if $\|f\|_1=1$.
Notice that if $f,g$ are probability measures then $f*g$
is also a probability measure.

\begin{prop}[{\cite[Theorem C]{Bo}}]\label{TheoremC} Assume $A\subset\Fp,
B\subset \Fps$ with $|A|\ge |B|$. For some $\cml>0$,
\be\label{bour1}
\sum_{b\in B} E(A, bA)\ll \( \min(p/|A|,|B| \)^{-\cml}|A|^3|B|.
\ee
\end{prop}

\begin{remark}
 An explicit version of Proposition \ref{TheoremC}, with $c_0=1/10430$, 
is given in \cite{BGlib}.
\end{remark}

Note that if $|A|<|B|$, we may decompose $B$ as a disjoint
union of at most $2|B|/|A|$ sets $B_j$ with $|A|/2< |B_j|\le|A|$ and apply
\eqref{bour1}
for each $B_j$. Hence
\[
\sum_{b\in B}E(A, bA) \ll
\Big[\min\Big(|A|,|B|, \frac p{|A|}\Big)\Big]^{-\cml} |A|^3|B|.
\]
Applying the \CS inequality we get
\be\label{bour3}
\sum_{b\in B} \Vert 1_A* 1_{bA}\Vert_2 \ll |A|^{3/2} \(
|A|^{-\cml/2}|B|+ |B|^{1-\cml/2} + p^{-\cml/2}|A|^{\cml/2}|B|\).
\ee

\begin{remark} It would be interesting to find best possible value for $\cml$
in Proposition~\ref{TheoremC}. The example $A=B=\{1,\dots,[\sqrt p]\}$
shows that $\cml<1$.
\end{remark}

\begin{cor}\label{measener} For any $A\subset\Fp$ and a probability measure
$\lambda$ we have
$$\sum_{b\in\Fps}\lambda(b)\|1_A*1_{bA}\|_2 \ll
\left(\|\lambda\|_2+|A|^{-1/2}+|A|^{1/2}p^{-1/2}\right)^{\cml}|A|^{3/2}.$$
\end{cor}

\begin{proof} Put $\lam(p)=0$, and let $b$ be a permutation of $\{1,\ldots,
p\}$ such that $\lam(b_1)\ge\dots\ge\lam(b_{p})=0$.
By \eqref{bour3}, for $1\le j\le p-1$ we have $S_j\ll G_j$, where
\[
 S_j = \sum_{h=1}^j \|1_A*1_{bA}\|_2, \quad G_j := |A|^{3/2}
\(|A|^{-\cml/2}j + |A|^{\cml/2}p^{-\cml/2}j+ j^{1-\cml/2} \).
\]
Applying summation by parts,
\begin{align*}
\sum_{b\in\Fps}\lam(b)\|1_A*1_{bA}\|_2 &= \sum_{j=1}^p \lam(b_j) \(S_j-
  S_{j-1}\) = \sum_{j=1}^{p-1} S_j \( \lam(b_j) - \lam(b_{j+1}) \) \\
&\ll \sum_{j=1}^{p-1} G_j \( \lam(b_j) - \lam(b_{j+1}) \) =
  \sum_{j=1}^{p-1} \lam(b_j) \(G_j-G_{j-1}\) \\
&= |A|^{3/2} \Big[ |A|^{-\cml/2}+p^{-\cml/2}|A|^{\cml/2} +
  O \Big( \sum_{j=1}^p \lam(b_j) j^{-\cml/2} \Big) \Big]
\end{align*}

Denote $u_0=\|\lambda\|_2^{-2}$. Notice that $1\le u_0\le p$
since $\|\lambda\|_1=1$.  Separately considering $j\le u_0$ and
$j>u_0$ and using the \CS inequality, we get
\[
 \sum_{j=1}^p \lam(b_j) j^{-\cml/2} \le \| \lam \|_2 \( \sum_{j\le u_0}
j^{-\cml}\)^{1/2} + u_0^{-\cml/2} =O \( \| \lam \|_2^{\cml} \).
\qedhere
\]
\end{proof}

Although Corollary~\ref{measener} suffices for the purposes of this paper,
a further generalization of Proposition~\ref{TheoremC}
might be useful.   For $z\in\Fps$ we define a function
$\rho_z[f]$ by $\rho_z[f](x)=f(x/z)$.

\begin{thm}\label{measener2} Let $\lambda,\mu$ be probability measures on
$\Fp$.
Then
$$
\sum_{b\in\Fps}\lambda(b)\|\mu*\rho_b[\mu]\|_2 \ll
\left(\|\lambda\|_2+\|\mu\|_2+\|\mu\|_2^{-1}p^{-1/2}\right)^{\cml/7}
\|\mu\|_2.
$$
\end{thm}

\begin{proof} Using a parameter $\Delta\ge1$ which will be specified later
we define the sets
$$
A_-=\{x:\,\mu(x)\ge\|\mu\|_2^2\Delta\},\quad
A_+=\{x:\,\mu(x)<\|\mu\|_2^2\Delta^{-2}\},\quad
A=\Fp\setminus A_-\setminus A_+.
$$
Decompose $\mu=\mu_-+\mu_0+\mu_+$ where
$$
\mu_-=\mu 1_{A_-},\quad \mu_0=\mu 1_A,\quad\mu_+=\mu 1_{A_+}.
$$
The contribution to the sum in the theorem from $\mu_-$ and $\mu_+$ is
negligible.  First,
\be\label{mu-}
\|\mu_-\|_1\le \frac{1}{\Delta \|\mu\|_2^2}
  \sum_{x\in A_-}\mu(x)^2 \le \Delta^{-1}.
\ee
and
\be\label{mu+}
\|\mu_+\|_2\le\|\mu\|_2\Delta^{-1}\|\mu_+\|_1^{1/2}\le\|\mu\|_2\Delta^{-1}.
\ee
Using Young's inequality (cf \cite{TaoVu}, Theorem~4.8), we find that
\be\label{contrmu-}\begin{split}
\sum_{b\in\Fps}\lambda(b)\|\mu_-*\rho_b[\mu]\|_2
&\le\sum_{b\in\Fps}\lambda(b)\|\mu_-\|_1\|\rho_b[\mu]\|_2 \\
&\le\sum_{b\in\Fps}\lambda(b)\Delta^{-1}\|\mu\|_2\le\Delta^{-1}\|\mu\|_2,
\end{split}\ee
\be\label{contrmu+}\begin{split}
\sum_{b\in\Fps}\lambda(b)\|\mu_+*\rho_b[\mu]\|_2
&\le\sum_{b\in\Fps}\lambda(b)\|\mu_+\|_2\|\rho_b[\mu]\|_1\\
&\le\sum_{b\in\Fps}\lambda(b)\Delta^{-1}\|\mu\|_2\le\Delta^{-1}\|\mu\|_2,
\end{split}\ee
Similarly,
\be\label{contrrho+}
\sum_{b\in\Fps}\lambda(b) \|\mu_0*\rho_b[(\mu_-+\mu_+)]\|_2
\le2\Delta^{-1}\|\mu\|_2,
\ee
So, it suffices to estimate the contribution of $\mu_0$.  We have
$$
1=\|\mu\|_1 \ge \sum_{x\in A}\mu(x)\ge|A|\|\mu\|_2^2\Delta^{-2}.
$$
Hence, $|A|\le\|\mu\|_2^{-2}\Delta^2$. Now we can use
Corollary~\ref{measener}:
\begin{align*}
\sum_{b\in\Fps}\lambda(b) & \|\mu_0*\rho_b[\mu_0]\|_2
\le\|\mu\|_2^4\Delta^2\sum_{b\in\Fps}\lambda(b)\|1_A*1_{bA}\|_2\\
&\ll \|\mu\|_2^4\Delta^2
\left(\|\lambda\|_2^{\cml}+|A|^{-\cml/2}+|A|^{\cml/2}p^{-\cml/2}\right)
|A|^{3/2}\\
&\le \|\mu\|_2^4\Delta^2\left(\|\lambda\|_2^{\cml}\|\mu\|_2^{-3}\Delta^3
+\|\mu\|_2^{-3+\cml}\Delta^{3-\cml}+\|\mu\|_2^{-3-\cml}\Delta^{3+\cml}\right)\\
&\le \Delta^6\|\mu\|_2
\left(\|\lambda\|_2^{\cml}
+\|\mu\|_2^{\cml}+\|\mu\|_2^{-\cml}p^{-\cml/2}\right).
\end{align*}
Combining the last inequality with \eqref{contrmu-} -- \eqref{contrrho+}
we get
$$\sum_{b\in\Fps}\lambda(b)\|\mu*\rho_b[\mu]\|_2\le
4\Delta^{-1}\|\mu\|_2+O(\Delta^6\|\mu\|_2S),$$
where
$$
S=\|\lambda\|_2^{\cml}+\|\mu\|_2^{\cml}+\|\mu\|_2^{-\cml}p^{-\cml/2}.
$$
Taking $\Delta=\max(1,S^{1/7})$ completes the proof of the theorem.
\end{proof}

%
%
\section{A sumset estimate in product sets}\label{Sec5}
%
%

The main result of this section is the following.
\begin{thm}\label{sumset} Let $r,M\in\NN, M\ge2$ and
$\cC=\cC_{M,r}=\{0,\dots,M-1\}^r$. Let $\tau=\tau_M$ be the solution
of the equation
$$\left(\frac1M\right)^{2\tau}+\left(\frac{M-1}M\right)^{\tau}=1.$$
Then for any subsets $A,B\subset\cC$ we have
\be\label{A+B}
|A+B|\ge(|A||B|)^\tau.
\ee
\end{thm}

Observe that for $A=B=\cC$ we have $|A+B|=|A|^{\tau'}|B|^{\tau'}$
where
$$\tau'=\tau'_M=\frac{\log(2M-1)}{2\log M}.$$
By Theorem \ref{sumset}, $\tau\le\tau'$. On the other hand, $\tau>1/2$.
If $M\to\infty$ then
\be\label{asymptau}
u^{2\tau}=1-(1-u)^\tau\sim\frac u2,\quad 2\tau-1\sim\frac{\log2}{\log M}
\sim2\tau'-1.
\ee
So, the asymptotic behavior of $2\tau_M-1$ as $M\to\infty$ is sharp.
Likely, inequality \eqref{A+B} holds with $\tau=\tau'$.  This was proved 
in the case $M=2$ by Woodall \cite{Wood}.

Results of a similar spirit, concerning addition of subsets of $\Fp^r$
and related groups, are considered in \cite{BL}.

For positive integers $K,L$ we define an $UR-$path as a sequence
of pairs of integers
$\CP=((i_1,j_1)=(0,0),\ldots,(i_{K+L-1},j_{K+L-1})=(K-1,L-1))$
such that for any $n$ either
$i_{n+1}=i_n+1, j_{n+1}=j_n$, or $i_{n+1}=i_n, j_{n+1}=j_n+1$.

\begin{lem}\label{ordered} Let $KL\le M^2$, $u_0\ge\dots\ge u_{K-1}\ge0$,
$v_0\ge\dots\ge v_{L-1}\ge0$, $\tau=\tau_M$.
Then there exists an $UR-$path $\CP$ such that
\be\label{ordineq}
\sum_{n=1}^{K+L-1}\left(u_{i_n}v_{j_n}\right)^\tau
\ge\left(\sum_{i=0}^{K-1}u_i\right)^\tau\left(\sum_{j=0}^{L-1}v_j\right)^\tau.
\ee
\end{lem}

\begin{proof} We proceed by induction on $K+L$. For $K=1$ or $L=1$
the assertion is obvious. We prove it for $K,L$ with
$\min(K,L)\ge2$, $KL\le M^2$ supposing that it holds for $(K,L)$
replaced by $(K-1,L)$ and $(K,L-1)$. Without loss of generality
we assume that
$$\sum_{i=0}^{K-1}u_i = \sum_{j=0}^{L-1}v_j =1.$$
By the induction supposition, there exists an $UR-$path $\CP$ such that
$i_1=1, j_1=0$ and
$$\sum_{n=2}^{K+L-1}\left(u_{i_n}v_{j_n}\right)^\tau
\ge\left(\sum_{i=1}^{K-1}u_i\right)^\tau\left(\sum_{j=0}^{L-1}v_j\right)^\tau
=(1-u_0)^\tau.$$
Therefore,
\[
S :=\max_{\CP}\sum_{n=1}^{K+L-1}\left(u_{i_n}v_{j_n}\right)^\tau
\ge(u_0v_0)^\tau+(1-u_0)^\tau.
\]
Similarly, $S \ge (u_0v_0)^\tau+(1-v_0)^\tau.$
Thus, $S \ge w^{2\tau}+(1-w)^\tau$
where
\[
 w=(u_0v_0)^{1/2}\ge(KL)^{-1/2}\ge 1/M.
\]
The function $f(x)=x^{2\tau}+(1-x)^\tau-1$ has negative third derivative
on $[0,1]$ and $f(0)=f(1/M)=f(1)=0$.  By Rolle's theorem, $f$ has no other
zeros on $[0,1]$, and since $f(u)>0$ for $u$ close to 1, $f(x)\ge0$
for $1/M\le x\le1$. Therefore, $f(w)\ge 0$ as desired.
\end{proof}

We will need Lemma~\ref{ordered} only for $K=L=M$
(although for the proof it was convenient to have varying $K,L$).

\begin{lem}\label{unordered} Let $U_0,\dots,U_{M-1}$, $V_0,\dots,V_{M-1}$
be non-negative numbers, and $\tau=\tau_M$. Then
\be\label{unordineq}
\sum_{\mu=0}^{2M-2}\max_{\substack{\kappa+\lambda=\mu,\\ \kappa\ge0,
\lambda\ge0}} (U_\kappa V_\lambda)^\tau
\ge\left(\sum_{\kappa=0}^{M-1}U_\kappa\right)^\tau
\left(\sum_{\lambda=0}^{M-1}V_\lambda\right)^\tau.
\ee
\end{lem}

Lemma~\ref{unordered} has some similarity
with inequality~(2.1) from \cite{Pre}.

\begin{proof} We order $U_0,\dots,U_{M-1}$ and $V_0,\dots,V_{M-1}$
in the descending order $u_0\ge\dots\ge u_{M-1}$ and
$v_0\ge\dots\ge v_{M-1}$, respectively, where for some permutations
$\pi$ and $\sigma$ of the set $\{0,\dots,M-1\}$ we have
$u_i=U_{\pi_i}, v_j=V_{\sigma_j}$. We consider an arbitrary
$UR-$path $\CP$ with $K=L=M$.  Since $|\{\pi_{i_1},\dots,\pi_{i_n}\}|=i_n+1$
and $|\{\sigma_{j_1},\dots,\sigma_{j_n}\}|=j_n+1$,
$$
|\{\pi_{i_1},\dots,\pi_{i_n}\}
+\{\sigma_{j_1},\dots,\sigma_{j_n}\}| \ge i_n+j_n+1.
$$
Consequently, there is a permutation $\psi$ of
$\{0,\dots,2M-2\}$ so that
$$
\psi(n-1)\in\{\pi_{i_1},\dots,\pi_{i_n}\}
+\{\sigma_{j_1},\dots,\sigma_{j_n}\} \qquad (1\le n\le 2M-1).
$$
Thus, for some $\kappa_0\in\{\pi_{i_1},\dots,\pi_{i_n}\}$ and
$\lambda_0\in\{\sigma_{j_1},\dots,\sigma_{j_n}\}$ we have
$$
\max_{\substack{\kappa+\lambda=\psi(n-1),\\ \kappa\ge0,
\lambda\ge0}} (U_\kappa V_\lambda)^\tau\ge (U_{\kappa_0} V_{\lambda_0})^\tau.
$$
But $U_{\kappa_0}=u_i$ for some $i\in\{i_1,\dots,i_n\}$.
Recalling that $i_1\le i_2\le\dots$ and $u_1\ge u_2\ge\dots$ we obtain
$U_{\kappa_0}\ge u_{i_n}$. Similarly, $V_{\lambda_0}\ge v_{j_n}$. Therefore,
$$
\max_{\substack{\kappa+\lambda=\psi(n-1),\\ \kappa\ge0,
\lambda\ge0}} (U_\kappa V_\lambda)^\tau\ge\left(u_{i_n}v_{j_n}\right)^\tau
$$
and
$$
\sum_{\mu=0}^{2M-2}\max_{\substack{\kappa+\lambda=\mu,\\ \kappa\ge0,
\lambda\ge0}} (U_\kappa V_\lambda)^\tau
=\sum_{n=1}^{2M-1}\max_{\substack{\kappa+\lambda=\psi(n-1),\\ \kappa\ge0,
\lambda\ge0}} (U_\kappa V_\lambda)^\tau
\ge\sum_{n=1}^{2M-1}\left(u_{i_n}v_{j_n}\right)^\tau,
$$
and the result follows from Lemma~\ref{ordered}.
\end{proof}

Now we are ready to prove Theorem~\ref{sumset}.
We proceed by induction on $r$. For $r=0$ the set $\cC_{M,r}$ is a
singleton, and there is nothing to prove. Now suppose that the assertion holds
for $r$ replaced by $r-1\ge0$. We consider arbitrary subsets
$A,B\subset\cC=\cC_{M,r}$. For $i=0,\dots,M-1$ we denote
$$A_i=\{(x_1,\dots,x_{r-1}):\,(x_1,\dots,x_{r-1},i)\in A\},$$
$$B_i=\{(x_1,\dots,x_{r-1}):\,(x_1,\dots,x_{r-1},i)\in B\}.$$
Let $D=A+B$. For $n=0,\dots,2M-2$ we denote
$$D_n=\{(x_1,\dots,x_{r-1}):\,(x_1,\dots,x_{r-1},n)\in D\}.$$
Observe that
$$|A|=\sum_i|A_i|,\quad B=\sum_j|B_j|,\quad D=\sum_n|D_n|.$$
For any $n=0,\dots,2M-2$ we have
$$|D_n|\ge\max_{\substack{i+j,\\ i\ge0, j\ge0}} |A_i+B_j|.$$
By the induction supposition, $|A_i+B_j|\ge(|A_i||B_j|)^\tau$. Hence,
$$|D_n|\ge\max_{\substack{i+j,\\ i\ge0, j\ge0}} (|A_i||B_j|)^\tau.$$
Applying Lemma~\ref{unordered},
\[
|D|=\sum_n|D_n|\ge\sum_n\max_{\substack{i+j,\\ i\ge0, j\ge0}}
(|A_i||B_j|)^\tau\ge\left(\sum_i|A_i|\right)^\tau
\left(\sum_j|B_j|\right)^\tau = (|A||B|)^\tau.
\qedhere
\]
The proof of Theorem \ref{sumset} is complete.

\begin{cor}\label{sumset1} Let $m$ be a positive integer.
For the set $\CB\subset\Fp$ defined in
Section~\ref{Sec2} and for any subset $B\subset\CB, |B|>p^{1/4}$ we have
$|B-B|\ge p^{\beta/5}|B|$.
\end{cor}

\begin{proof}
The set $-B$ is a translate of some set $B'\subset\CB$, and
$\CB$ is Freiman isomorphic to $\cC_{M,r}$. Hence, for
any $B\subset\CB$ we have $|B-B|=|B+B'|\ge|B|^{2\tau_M}$. If $|B|>p^{1/4}$
then $|B-B|\ge|p|^{(2\tau_M-1)/4}|B|$. By \eqref{asymptau} and a short calculation using $M\ge 2^{15}$, $p^{(2\tau_M-1)/4}\ge p^{\beta/5}$.
\end{proof}

\begin{cor}\label{sumset2} Fix $m\in\NN$ and let $p\ge p(m)$ be a sufficiently
large prime. Let $\CB\subset\Fp$ be the set defined in
Section~\ref{Sec2}. Then for any subset $S\subset\CB, |S|>p^{1/3}$ we have
$E(S,S)\le p^{-\beta/50}|S|^3$.
\end{cor}

\begin{proof} Let $E(S,S)=|S|^3/K$. By Corollary~\ref{BSG2}, there is a set
$B\subset S$ such that $|B|\ge|S|/(20K)$ and $|B-B|\le10^7K^9|S|$.
If $K\le p^{\beta/50}<p^{1/24}$ and $p$ is so large that
$10^7\le p^{\beta/50}$ then we get contradiction with
Corollary~\ref{sumset1}.
\end{proof}

%
%
\section{The proof of Lemma~\ref{main}}\label{Sec7}
%
%

We may assume $\eps_1>0$, otherwise
there is nothing to prove.  Adopt the notation ($A_i, M_i, \Omega_i(a)$)
from Section \ref{Sec2}.
If $|A_1|M_1< p^{1/2-\gamma/10}$, then by \eqref{balanceM},
$|S(A_1,A_2)|\le2p^{1-\gamma/10}$ and \eqref{balanced} holds (recall
that $\cml<1$, hence $\eps_1<\g/20$).
Thus, we can assume that
$|A_1|M_1\ge p^{1/2-\gamma/10}$, which implies, by \eqref{Aupper}, that
\be\label{lowernewM1}
M_1\ge p^{1/2-\alpha-\gamma/10}.
\ee

\begin{lem}\label{usingener} For any $\theta\in\Fp^*$,
$B_1\subset\Fp$, $B_2\subset\Fp$ we have
\[
\Bigg|\sum_{\substack{b_1\in B_1\\b_2\in B_2}}
e_p\left(\theta(b_1-b_2)^2\right)\Bigg|
\le|B_1|^{1/2}E(B_1,B_1)^{1/8} |B_2|^{1/2}E(B_2,B_2)^{1/8}p^{1/8}.
\]
\end{lem}

\begin{proof} Let $W$ denote the double sum over $b_1,b_2$.  By the \CS inequality,
\begin{align*}
|W|^2 &\le |B_1| \sum_{b_1\in B_1} \Bigg| \sum_{b_2\in B_2}
e_p\left(\theta(b_1-b_2)^2\right)\Bigg|^2 \\
&=|B_1| \sum_{b_2,b_2'\in B_2}\sum_{b_1\in B_1}
e_p\left(\theta\left(b_2^2-(b_2')^2-2b_1(b_2-b_2')\)\).
\end{align*}
Another application of the \CS inequality gives
\begin{align*}
|W|^4 &\le |B_1|^2|B_2|^2  \sum_{b_2,b_2'\in B_2}\Bigg| \sum_{b_1}
e_p \( 2\theta b_1 (b_2-b_2') \) \Bigg|^2 \\
&= |B_1|^2|B_2|^2 \sum_{x,y\in\Fp}\lambda_x\mu_y e_p(-2\theta xy),
\end{align*}
where
$$
\lambda_x=1_{B_1} * 1_{(-B_1)}(x),\qquad
\mu_y=1_{B_2} * 1_{(-B_2)}(y).
$$
A third application of the \CS inequality, followed by Parseval's identity
yields a well-known inequality (cf. \cite{Vin}, Problem~14(a) for Chapter~6)
\begin{align*}
 \Bigg|\sum_{x,y\in\Fp}\lambda_x\mu_y e_p(-2\theta xy)\Bigg|^2 &\le
\| \lambda \|_2^2 \sum_{x\in\Fp} \Bigg| \sum_{y\in\Fp} \mu_y e_p(-2\theta xy)
\Bigg|^2 \\
&=p \| \lambda \|_2^2 \| \mu \|_2^2 = p E(B_1,B_1) E(B_2,B_2).
\qedhere
\end{align*}
\end{proof}

By \eqref{lowernewM1}, $|\Omega_i(a_i)| \ge p^{1/3}$, and by
Lemma~\ref{usingener} and ~\eqref{Benergy},
$$
\Bigg|\sum_{\substack{b_1\in\Omega_1(a_1)\\b_2\in\Omega_2(a_2)}}
e_p\left(\frac{(b_1-b_2)^2}{4(a_1-a_2)}\right)\Bigg|
\le|\Omega_1(a_1)|^{7/8}|\Omega_2(a_2)|^{7/8}p^{1/8-\gamma/4}.
$$
Next, by \eqref{balanceM}, we have
$$
|S(A_1,A_2)|\leq4|A_1|^{1/8}|A_2|^{1/8}p^{1-\gamma/4}.
$$
Thus, if $|A_1|<p^{\gamma/2}$ and $|A_2|<p^{\gamma/2}$, then
$|S(A_1,A_2)|\le 4p^{1-\g/8}$ and
\eqref{balanced} follows.  Otherwise,
without loss of generality we may assume that
\be\label{A2large}
|A_2|\ge p^{\gamma/2}.
\ee
The following lemma gives the necessary estimates to complete the
proof of Lemma \ref{main}.
For $a_1\in A_1$, set
\[
 T(A,B)=T_{a_1}(A,B)=\sum_{\substack{b_1\in B\\a_2\in A, b_2\in\Omega_2(a_2)}}
\left(\frac{a_1-a_2}p\right)
e_p\left(\frac{(b_1-b_2)^2}{4(a_1-a_2)}\right)
\]


\begin{lem}\label{2ndcase}
If $a_1\in A_1$, $0 < \gamma \le \min(\a,\frac{1}{3m})$,
conditions \eqref{balanceM} and \eqref{A2large}
are satisfied and a set $B\subset\Fp$ is such that
\be\label{Blarge}
p^{1/2-6\alpha}\le|B|\le p^{1/2}
\ee
and
\be\label{B-Bsmall}
|B-B|\le p^{28\alpha}|B|,
\ee
then
\be\label{crucial}
|T(A_2,B)| \le |B|p^{(1/2)-\eps_2},\qquad \eps_2=\frac{\cml \g}{20}
-\frac{42\alpha}{m}.
\ee
\end{lem}

\begin{remark}
 The proof of Lemma \ref{2ndcase} applies to more general sums, e.g. in
 $T(A,B)$ one may replace
the Legendre symbol $(\frac{a_1-a_2}{p})$ with arbitrary complex
numbers $\psi(a_1,a_2)$ with modulus $\le 1$, and one may replace
$\frac{1}{a_1-a_2}$ with different quantities $g(a_1,a_2)$ having
the dissociative property (the analog of \eqref{Asparse} holds).
\end{remark}

Postponing the proof of Lemma~\ref{2ndcase}, we show first how to deduce
Lemma \ref{main}.

We take a maximal subset $B_0\subset\Omega_1(a_1)$ so that
\eqref{crucial} holds for $B=B_0$. Denote $B_1=\Omega_1(a_1)\setminus B_0$.
By Lemma~\ref{usingener}, \eqref{balanceM}, and  \eqref{Aupper} we have
\begin{align*}
|T_{a_1}(A_2,B_1)| &\le
 \sum_{a_2\in A_2}|B_1|^{1/2}E(B_1,B_1)^{1/8}
 |\Omega_2(a_2)|^{1/2}E(\Omega_2(a_2),\Omega_2(a_2))^{1/8}p^{1/8}\\
&\le |A_2|\, |B_1|^{1/2}E(B_1,B_1)^{1/8} M_2^{7/8}p^{1/8}\\
&\le 2 |B_1|^{1/2}E(B_1,B_1)^{1/8}p^{(9/16)+(\alpha/8)}.
\end{align*}
Consider the case when
\be\label{EBcase1}
E(B_1,B_1)\le p^{-3\alpha}M_1^3.
\ee
Then we have, due to \eqref{balanceM},
\be\label{case1est}
|T_{a_1}(A_2,B_1)|
\le2 M_1^{7/8}p^{(9/16)-\alpha/4}.
\ee
Now assume that \eqref{EBcase1} does not hold.  By
\eqref{balanceM}, we get
$$
|B_1|>p^{-\alpha}M_1,\quad E(B_1,B_1)\ge p^{-3\alpha}|B_1|^3.
$$
Applying now Corollary~\ref{BSG2} and \eqref{balanceM} we obtain the existence of a set
$B_1'\subset B_1$ such that
$$
|B_1'|\ge \frac{M_1}{20p^{4\alpha}} \ge \frac{p^{1/2-5\a-\g/10}}{20}
\ge p^{1/2-6\a}
$$
and $|B_1'-B_1'|\le 10^7p^{27\alpha}|B_1| \le p^{28\a} |B_1|$.
Using Lemma~\ref{2ndcase} we get inequality~\eqref{crucial} for $B=B_1'$.
Therefore, \eqref{crucial} is also satisfied
for $B=B_0\cup B_1'$, contradicting the choice of $B_0$.

Thus, we have shown that \eqref{EBcase1}
must hold. Using \eqref{crucial} for $B=B_0$ and \eqref{case1est} we get
$$
|T_{a_1}(A_2,\Omega_1(a_1))|
\le M_1p^{(1/2)-\eps_2} + 2M_1^{7/8}p^{(9/16)-\alpha/4}.
$$
Summing on $a_1\in A_1$ and using \eqref{Aupper} and \eqref{balanceM},
we obtain
\begin{align*}
|S(A_1,A_2)| &\le
|A_1|\left(M_1p^{(1/2)-\eps_2} + 2M_1^{7/8}p^{(9/16)-\alpha/4}\right)\\
&\le2p^{1-\eps_2}+4|A_1|^{1/8}p^{1-\alpha/4}\le2p^{1-\eps_2}+4p^{1-\alpha/8},
\end{align*}
completing the proof of Lemma \ref{main}.


\begin{proof}[Proof of Lemma \ref{2ndcase}]
By the \CS inequality we have
\[
|T(A_2,B)|^2\le\sqrt p\sum_{b_1,b\in B}|F(b,b_1)|,
\]
where
$$
F(b,b_1) = \sum_{\substack{a_2\in A_2\\b_2\in\Omega_2(a_2)}}
e_p\left(\frac{b_1^2-b^2}{4(a_1-a_2)}-\frac{b_2(b_1-b)}{2(a_1-a_2)}\right).
$$
Consequently, by H\"older's inequality,
\be\label{5.3}
|T(A_2,B)|^2\le\sqrt p|B|^{2-2/m}\left(\sum_{b_1,b\in B}|F(b,b_1)|^m
\right)^{\frac1m}.
\ee
Next,
\begin{align*}
\sum_{b_1,b\in B}|F(b,b_1)|^m &\le
\sum_{\substack{x\in B+B,\\y\in B-B}}\Bigg|
\sum_{\substack{a_2\in A_2,\\b_2\in\Omega_2(a_2)}}
e_p\left(\frac{xy}{4(a_1-a_2)}-\frac{b_2y}{2(a_1-a_2)}\right)\Bigg|^m\\
&\le \sum_{y\in B-B} \sum_{\substack{a_2^{(i)}\in A_2 \\ b_2^{(i)}\in\Omega_2(a_2^{(i)}) \\ 1\le i\le m}} \Bigg| \sum_{x\in B+B}
e_p\Bigg( \frac{xy}{4} \sum_{i=1}^{m/2} \Bigg[ \frac{1}{a_1-a_2^{(i)}}-
\frac{1}{a_1-a_2^{(i+m/2)}} \Bigg] \Bigg) \Bigg|.
\end{align*}
Hence, for some complex numbers $\eps_{y,\xi}$ of modulus $\le 1$,
\be\label{5.4}
\sum_{b_1,b\in B}|F(b,b_1)|^m \le M_2^m \sum_{y\in B-B} \; \sum_{\xi\in \Fp}
\lam(\xi) \eps_{y,\xi}  \sum_{x\in B+B} e_p(xy\xi/4),
\ee
where
\begin{align*}
\lambda(\xi) =\Bigg|\Bigg\{a^{(1)},\ldots,a^{(m)} \in A_2 : \;
\sum_{i=1}^{m/2} \( \frac{1}{a_1-a^{(i)}} - \frac{1}{a_1-a^{(i+m/2)}} \)=\xi
\Bigg\}\Bigg|.
\end{align*}

By \eqref{Asparse},
\be\label{lambda0}
\lambda(0)\le(m/2)!|A_2|^{m/2}.
\ee
Let
\[
 \zeta'(z)=\sum_{\substack{y\in B-B\\ \xi\in \Fps \\ y\xi=z}} \eps_{y,\xi}
\lam(\xi),
\qquad \zeta(z)=\sum_{\substack{y\in B-B\\ \xi\in \Fps \\ y\xi=z}} \lam(\xi).
\]
Then $|\zeta'(z)| \le \zeta(z)$.  By H\"older's inequality,
\be\label{Holder}\begin{split}
 \Bigg| \sum_{y\in B-B} \sum_{\xi\in \Fps}
\lam(\xi) \eps_{y,\xi}  &\sum_{x\in B+B} e_p(xy\xi/4) \Bigg| =
\Bigg| \sum_{\substack{x\in B+B \\ z\in \Fp}} \zeta'(z) e_p(xz/4) \Bigg| \\
&\le |B+B|^{3/4} \Bigg( \sum_{x\in\Fp} \Bigg| \sum_{z\in \Fp}
\zeta'(z) e_p(xz/4) \Bigg|^4 \Bigg)^{1/4} \\
&= |B+B|^{3/4} \Bigg( \sum_{x\in\Fp}  \Bigg| \sum_{z'\in \Fp}
(\zeta' * \zeta')(z')e_p(xz'/4) \Bigg|^2 \Bigg)^{1/4} \\
&= |B+B|^{3/4} \| \zeta' * \zeta'\|_2^{1/2} p^{1/4} \\
&\le |B+B|^{3/4} \| \zeta * \zeta \|_2^{1/2} p^{1/4}.
\end{split}\ee
As $\zeta(z)=\sum_{\xi}1_{B-B}(z/\xi)$, we have by the triangle inequality,
\be\label{zz}
\begin{split}
 \| \zeta * \zeta \|_2 &\le  \sum_{\xi,\xi'\in \Fps} \lam(\xi)\lam(\xi') \| 1_{\xi(B-B)} *
1_{\xi'(B-B)}\|_2\\
&=\sum_{\xi,\xi'\in \Fps} \lam(\xi)\lam(\xi') \| 1_{B-B} *
1_{(\xi'/\xi)(B-B)}\|_2.
\end{split}\ee

Define the probability measure $\lam_1$ by
\[
 \lam_1(\xi) = \frac{\lam(\xi)}{\|\lam\|_1} = \frac{\lam(\xi)}{|A_2|^m}.
\]
The sum $\sum_{\xi\in\Fp}\lambda(\xi)^2$ is equal to the number of solutions
of the equation
\[
\frac1{a_1-a^{(1)}}+\dots+\frac1{a_1-a^{(m)}}
-\frac1{a_1-a^{(m+1)}}-\frac1{a_1-a^{(2m)}}=0
\]
with $a^{(1)},\dots,a^{(2m)}\in A_2$. By \eqref{Asparse}, this
has only trivial solutions and thus
\be\label{sumlamsq}
\sum_{\xi\in\Fp}\lambda(\xi)^2\le m!|A_2|^m.
\ee

Now we are in position to apply Corollary~\ref{measener} which gives
for any $\xi'\in\Fps$
\begin{multline}\label{Cor5}
\sum_{\xi\in\Fps}\lambda_1(\xi) \|1_{B-B}*1_{(\xi'/\xi)(B-B)}\|_2\\
\ll \left(\|\lambda_1\|_2+|B-B|^{-1/2}
+|B-B|^{1/2}p^{-1/2}\right)^{\cml}|B-B|^{3/2}.
\end{multline}
By \eqref{A2large} and \eqref{sumlamsq},
$$
\|\lambda_1\|_2\le \sqrt{m!}p^{-m\gamma/4}.
$$
By \eqref{Blarge} and $\a<0.01$,
$$
|B-B| \ge |B| \ge p^{1/2-6\a} \ge p^{0.44}.
$$
On the other hand, it follows from \eqref{Blarge} and \eqref{B-Bsmall}
that
$$|B-B|\le p^{1/2+28\alpha} \le p^{0.78}.$$
Since $m\gamma \le 1/3$ we get
$$
\|\lambda_1\|_2+|B-B|^{-1/2}
+|B-B|^{1/2}p^{-1/2}\le\sqrt{m!}p^{-m\gamma/4}+p^{-0.1}\le p^{-m\gamma/5}.
$$
So, by \eqref{zz} and \eqref{Cor5},
\begin{align*}
\|\zeta*\zeta\|_2 &\le |A_2|^{2m} \sum_{\xi'\in\Fps} \lam_1(\xi')
\sum_{\xi\in\Fps}\lambda_1(\xi) \|1_{B-B}*1_{(\xi'/\xi)(B-B)}\|_2\\
&\ll |A_2|^{2m} p^{-(\cml/5)m\gamma}|B-B|^{3/2}.
\end{align*}
Subsequent application of \eqref{5.4}, \eqref{lambda0} and
\eqref{Holder}  gives
\begin{multline*}
 \sum_{b_1,b\in B}|F(b,b_1)|^m
\le (\tfrac{m}{2})!(M_2|A_2|)^m|A_2|^{-m/2}|B-B||B+B|\\
+O(M_2^m|A_2|^m|B-B|^{3/4}|B+B|^{3/4}p^{-(\cml/10)m\gamma} p^{1/4}).
\end{multline*}
Due to Lemma~\ref{PlunRuz}, condition~\eqref{B-Bsmall} implies
\[
|B+B|\le p^{56\alpha}|B|.
\]
By \eqref{Blarge}, $p^{1/4}\le |B|^{1/2} p^{3\a}$.
Recalling $\g\le \a$, \eqref{balanceM}, \eqref{A2large}
and \eqref{B-Bsmall}, we conclude that
\begin{align*}
\sum_{b_1,b\in B} |F(b,b_1)|^m
&\ll (\tfrac{m}{2})!(2\sqrt p)^mp^{-m\gamma/4}p^{84\alpha}|B|^2
+(2\sqrt p)^mp^{63\alpha}|B|^{3/2}p^{-(\cml/10)m\gamma}p^{1/4}\\
&\le |B|^2 p^{m/2-(\cml/10)m\g+84\a}.
\end{align*}
Plugging the last estimate into \eqref{5.3}, we get
\[
|T(A_2,B)|^2 \le \sqrt p|B|^{2-2/m}\left( |B|^2
p^{m/2-(\cml/10)m\gamma+84\alpha}\right)^{\frac1m}\\
\le|B|^2p^{1+84\alpha/m-(\cml/10)\gamma}.
\qedhere
\]
\end{proof}

%
%
\section{Thin sets with small Fourier coefficients}\label{Sec8}
%
%

Denote by $(a^{-1})_m$ the inverse of $a$ modulo $m$.  It is easy to
see for relatively prime integers $a,b$ that

\be\label{reciprocity}
\frac{(a^{-1})_b}{b} + \frac{(b^{-1})_a}{a} - \frac{1}{ab} \in \mathbb{Z}.
\ee

\begin{lem}\label{iter1}
Let $P\ge 4$, $S\ge 2$, and $R$ be a positive integer. Suppose
that for every prime $p\le P$, $S_p$ is a set of
integers in $(-p/2, p/2)$.
Suppose $q$ is a prime satisfying $q \ge R P^2$.  Then
the numbers
$r + s^{(p)} (p^{-1})_q$, where $1\le r \le R, P/2 < p \le P,
s^{(p)}\in S_p$, are distinct modulo $q$.
\end{lem}

\begin{proof} Suppose that
$$
r_1 + s_1^{(p_1)} (p_1^{-1})_q \equiv r_2 + s_2^{(p_2)} (p_2^{-1})_q \pmod{q}.
$$
Multiplying both sides by $p_1 p_2$ gives
$$
r_1 p_1 p_2 + p_2  s_1^{(p_1)} \equiv r_2 p_1 p_2 + p_1 s_2^{(p_2)} \pmod{q}.
$$
By hypothesis,
$$
\left| (r_1-r_2) p_1 p_2 +  p_2  s_1^{(p_1)} -  p_1 s_2^{(p_2)}
\right| < (R-1)P^2 + P^2 \le q,
$$
thus
$$
(r_1-r_2) p_1 p_2  = - p_2 s_1^{(p_1)} + p_1 s_2^{(p_2)}.
$$
The right side is divisible by $p_1 p_2$ and the
absolute value of the right side is $< p_1 p_2$, hence both sides are
zero, $r_1=r_2$, $p_1=p_2$ and $s_1^{(p_1)} = s_2^{(p_2)}$.
\end{proof}

For brevity, we write $e(z)$ for $e^{2\pi i z}$ is what follows.

\begin{lem}\label{iter2}
Let $P\ge 4$, $S\ge 2$, and $R$ be a positive integer. Suppose
that for every prime $p\in(P/2,P]$, $S_p$ is a multiset of
integers in $(-p/2,p/2)$, $|S_p|=S$ and $|f_{S_p}| \le \eps$.
Suppose $q$ is a prime satisfying $q > P$.  Then the multiset
$$
T = \{ r + s^{(p)} (p^{-1})_q : 1\le r \le R, P/2 < p \le P,
s^{(p)}\in S_p \}
$$
of residues modulo $q$, satisfies
\be\label{fT}
|f_T| \le \eps + \frac{2/\sqrt{3}}{R} + \frac{\log (q/3)}{V \log (P/2)},
\ee
where $V$ is the number of primes in $(P/2,P]$.
\end{lem}

\begin{proof}
Since $|f_T(k)| = |f_T(q-k)|$, we may assume without
  loss of generality that $1\le k < q/2$.  We have
$$
f_T(k) = A(k) \sum_{P/2<p \le P} B(p,k),
$$
where
$$
A(k) = \sum_{r\le R} e\pfrac{kr}{q}, \qquad B(p,k) =
\sum_{s\in S_p} e\pfrac{ks(p^{-1})_q}{q}.
$$
Trivially,
\be\label{trivial2}
|A(k)| \le \min\( R, \frac{2}{|e(k/q)-1|} \).
\ee
If $k\ge q/3$, we use the trivial bound $|B(p,k)| \le S$ and
conclude
$$
\frac{|f_T(k)|}{|T|} \le \frac{2}{R|e(k/q)-1|} \le \frac{2}{R|e(1/3)-1|} =
\frac{2/\sqrt{3}}{R}.
$$
Now assume $k\le q/3$.
If $p|k$, then $|B(p,k)| \le S$.
When $p \nmid k$, by \eqref{reciprocity},
\begin{align*}
|B(p,k)| &= \left| \sum_{s\in S_p} e\( - \frac{s k (q^{-1})_p}{p} +
\frac{ks}{pq} \) \right| \\
&\le |S_p| \max_{s\in S_p} \left| e\pfrac{ks}{pq} - 1 \right| + \left|
\sum_{s\in S_p} e\pfrac{sk(q^{-1})_p}{p} \right| \\
&\le \( \eps + |e(k/2q)-1| \) S.
\end{align*}
Since there are $\le \frac{\log k}{\log (P/2)}$ primes $p|k$ with
$p>P/2$, we have
\[
\sum_{P/2 < p \le P} |B(p,k)| \le  \( \eps +  |e(k/2q)-1| \)
SV + \frac{\log (q/3)}{\log (P/2)} S.
\]
Combining our estimates for $|A(k)|$ and $|B(p,k)|$, we arrive at
\begin{align*}
\frac{|f_T(k)|}{|T|} &\le \eps + \frac{\log(q/3)}{V \log(P/2)} +
\frac{2}{R} \left| \frac{e(k/2q)-1}{e(k/q)-1} \right| \\
&\le \eps + \frac{\log(q/3)}{V \log(P/2)} + \frac{2/\sqrt{3}}{R}.
\qedhere
\end{align*}
\end{proof}

For a specific choice of $S_p$, the inequality \eqref{fT} can be strengthened.

\begin{lem}\label{iter3}
Let $P\ge 4$ and $R$ be a positive integer. For every prime $p\in(P/2,P]$
denote by $S_p$ the set of all integers in $(-p/2, p/2)$.
Suppose $q$ is a prime satisfying $q > P$.  Then the multiset
$$
T = \{ r + s^{(p)} (p^{-1})_q : 1\le r \le R, P/2 < p \le P,
s^{(p)}\in S_p \}
$$
of residues modulo $q$ satisfies
\be\label{fT3}
|f_T| \le \frac{W}{2V} + \frac{W}{RV}
\left(1 + \frac{\log\left(1+\frac VW\right)}2\right).
\ee
where $V$ is the number of primes in $(P/2,P]$ and
$W= 4\frac{\log(q/2)}{\log(P/2)}$.
\end{lem}

\begin{proof} Again, we may assume without loss of generality that
$1\le k < q/2$. We use notation from the proof of Lemma \ref{iter2}.
If $p|k$, we use the trivial estimate $|B(p,k)| \le |S_p| \le P$.
Now there are $\le \frac{\log (q/2)}{\log (P/2)}$ primes $p|k$ with
$p>P/2$. When $p \nmid k$, by \eqref{reciprocity},
\begin{align*}
|B(p,k)|
&\le \left|\sum_{s= (1-p)/2}^{(p-1)/2}e\( - \frac{s k (q^{-1})_p}{p} +
\frac{ks}{pq} \) \right| 
=\frac{\left|e\( \frac{k}{q} \)-1\right|}
{\left|e\( - \frac{k (q^{-1})_p}{p} +\frac{k}{pq} \)-1 \right|} \\
&\le\frac{\left|e\( \frac{k}{q} \) -1\right|}
{\left|e\( - \frac{2 |k(q^{-1})_p|-1}{2p} \)-1 \right|}
\le\frac{\left|e\( \frac{k}{q} \) -1\right|}
{\left|e\( - \frac{2 |k(q^{-1})_p|-1}{2P} \)-1 \right|},
\end{align*}
where it is assumed that $k(q^{-1})_p\in(-p/2,p/2)$.
For $a=1,\dots,[(P-1)/2]$ we denote
$$P_a=\{p\in(P/2,P]: |k(q^{-1})_p|=a\}.$$
Taking into account that $|e(u)-1|^{-1}\le1/(4u)$ for $u\in(0,1/2]$ we get
\be\label{Bsumest}
\sum_{p \nmid k}|B(p,k)| \le 
\frac{P}2\left|e\( \frac{k}{q} \) -1\right|\sum_a|P_a|\frac1{2a-1}.
\ee

If $k(q^{-1})_p=\pm a$ then $k\pm aq$ is divisible by $p$. But
$|k\pm aq| \le Pq/2$.
Therefore, the number of prime divisors $p>P/2$ of any number $k\pm aq$
is at most $\frac{\log q}{\log P/2}+1$ and for any $a$ we get
$$
|P_a| \le  2\left[\frac{\log q}{\log(P/2)}\right]+2 \le W.
$$
Let $A=[V/W]+1$. We have
\begin{align*}
\sum_a|P_a|\frac1{2a-1} &\le\sum_{a\le A}|P_a|\frac1{2a-1}
+\left(V-\sum_{a\leq A}|P_a|\right)\frac1{2A+1}\\
&\le\sum_{a\le A}W\frac1{2a-1}
+\left(V-\sum_{a\leq A}W\right)\frac1{2A+1}\le\sum_{a\le A}W\frac1{2a-1}\\
&\le W\left(1+\frac{\log A}2\right)
\le W\left(1+\frac{\log\left(1+\frac{V}{W}\right)}2\right).
\end{align*}
Combining our estimates for $|A(k)|$ and $|B(p,k)|$ (\eqref{trivial2} and
 \eqref{Bsumest}), we arrive at
\begin{align*}
\frac{|f_T(k)|}{|T|} &\le \frac{2\log(q/2)}{V \log(P/2)}
+ \frac{PW/2}{R(P-2)V/2}\left(1+\frac{\log\left(1+\frac{V}{W}\right)}2\right)\\
&= \frac{W}{2V}
+ \frac{W}{RV}\left(1+\frac{\log\left(1+\frac{V}{W}\right)}2\right).
\qedhere
\end{align*}
\end{proof}

\begin{remark} Applying Lemma \ref{iter2} for all primes $q$ in a dyadic
interval, we can then feed these multisets $T=T_q$
back into the lemma and iterate.
\end{remark}

Using explicit estimates for counts of prime numbers \cite{RS}, we have

\begin{prop}\label{primes}
For $P\ge 250$, there are more than $\frac{2P}{5\log (P/2)}$
primes in $(P/2,P]$.  For any $P>2$, there are at most $0.76 P/\log P$
  primes in $(P/2,P]$.
\end{prop}

Using Proposition \ref{primes} we obtain a more convenient version of
Lemma \ref{iter3}.

\begin{lem}\label{iter4} Let $P\ge 250$. For every prime $p\in(P/2,P]$
denote by $S_p$ the set of all nonzero integers in $(-p/2, p/2)$.
Suppose $q$ is a prime satisfying $q > P$ and suppose
$R\ge 1+\log(1+0.26P/\log(2q))/2$
is a positive integer. Then the multiset
$$
T = \{ r + s^{(p)} (p^{-1})_q : 1\le r \le R, P/2 < p \le P,
s^{(p)}\in S_p \}
$$
of residues modulo $q$ satisfies
\be\label{fT4}
|f_T| \le 15\frac{\log q}P.
\ee
\end{lem}

\begin{proof} We use the notation of Lemma \ref{iter3}. By Proposition
\ref{primes} we have
\be\label{estWV}
\frac{W}{2V}\le 5\frac{\log q}P.
\ee
On the other hand, using Proposition \ref{primes} again we get
$$
\frac VW\le\frac{0.76P/\log P}{4\log(q/2)/\log(P/2)} \le 0.19
 \frac{P}{\log(q/2)} \le 0.26 \frac{P}{\log (2q)}.
$$
Hence,
$$
R \ge 1 + \frac{\log\left(1+\frac VW\right)}2.
$$
Now the inequality \eqref{fT4} follows from \eqref{estWV} and \eqref{fT3}.
\end{proof}

Using just one iteration one can get the following effective result
on thin sets with small Fourier coefficients, of nearly the same strength
as \eqref{Katz}.

\begin{cor}\label{oneiter}
For sufficiently large prime $N$ and $\mu$ such that
$N^{-1/2}\log^2N \le \mu <1$
there is a set $T$ of residues modulo $N$ so that
$$
|f_T| \le \mu, \qquad
|T| =O\( \frac{L_1^2}{\mu^2}\pfrac{1+\log(1/\mu)}{L_2+\log(1/\mu)} \).
$$
\end{cor}

\begin{proof}
We choose $P=(15/\mu)\log N$ and
$$
R = \left[2 + \frac{\log\left(1+5/\mu\right)}2\right]
\ge 1 + \frac{\log\left(1+\frac{0.26P}{\log N}\right)}2.
$$
Clearly, $R\ll 1+\log(1/\mu)$.
Let $T$ be the multiset constructed in Lemma \ref{iter4}.
We have $|f_T| \le \mu$. By Lemma \ref{iter1}, $T$ is a set. Moreover,
\[
|T|\ll P^2\frac{1+\log(1/\mu)}{\log P} \ll \frac{P^2(1+\log(1/\mu))}
{L_2+\log(1/\mu)}.\qedhere
\]
\end{proof}

\begin{proof}[Proof of Theorem \ref{thm1}]
We choose real parameters $P_0$, $P_1$ and positive integers $R_0$, $R_1$ so
that
\be\label{P1P2}
P_0 \ge 250, \qquad P_1 \ge 2 R_0 P_0^2, \qquad N \ge R_1 P_1^2, \qquad
R_0 \ge 1 + \frac{\log\left(1+\frac{0.26P_0}{\log P_1}\right)}2
\ee
and also
\be\label{errors}
\frac{2/\sqrt{3}}{R_1} + 15\frac{\log P_1}{P_0} +  \frac{5\log N}{2P_1}
\le \mu.
\ee
For $P_0/2 < p \le P_0$, let $S_p$ be the set of integers in $(-p/2,p/2)$. 
By Lemmas  \ref{iter1},
\ref{iter4} and \eqref{P1P2}, for each prime $q\in (P_1/2,P_1]$, there is
a set $T=S_q$ of residues modulo $q$ such that
$$
|f_{S_q}| \le 15 \frac{\log (P_1)}{P_0}=:\eps_1.
$$
By an  application of Lemmas  \ref{iter1} and \ref{iter2} with $P=P_1$,
$\eps=\eps_1$, $q=N$, and $S=R_0\sum_{P_0/2<p\le P_0} p$, together with
\eqref{errors}, there is a set $T$ of residues modulo $N$ so that
$$
|f_T| \le \eps_1 + \frac{2/\sqrt{3}}{R_1} +  \frac{5\log N}{2P_1} \le \mu.
$$
Using Proposition \ref{primes}, we find that
$$
|T| \le (0.76)^2 R_0 R_1 \frac{P_1 P_0^2}{(\log P_0)(\log P_1)}.
$$
Recalling that $1/\mu\in \NN$, we now take
\begin{align*}
R_0 = \left[ 2+\log(1+13/\mu)/2\right], \qquad R_1 = 4/\mu,\\
 P_1 = (8/\mu)\log N, \qquad P_0=(45/\mu)\log P_1
\end{align*}
so that \eqref{errors} follows immediately.  The condition \eqref{KN}
implies \eqref{P1P2} for large enough $N$.
\end{proof}

\begin{remark} Theorem \ref{thm1} supersedes Corollary \ref{oneiter}
for $\mu\gg L_1^{-1/2}L_2^{1/2}$.
\end{remark}

%
%
\section{An explicit construction for Tur\'an's problem}\label{Sec9}
%
%

\begin{proof}[Proof of Theorem \ref{thm2}]
We follow the proof of Theorem \ref{thm1} and Lemma \ref{iter2}.
We choose real parameters $P_0$, $P_1$ and a positive integer $R_0$, so
that
\be\label{P1P2-2}
P_0 \ge 250, \qquad P_1 > 2 P_0^2, \qquad
R_0 \ge 1 + \frac{\log\left(1+\frac{0.26P_0}{\log P_1}\right)}2
\ee
and also
\be\label{errors2}
15\frac{\log P_1}{P_0} + \frac{5\log N}{2P_1} \le \mu.
\ee
For $P_0/2 < p \le P_0$, let $S_p$ be the set of integers in $(-p/2,p/2)$.
 By Lemma
\ref{iter4} and \eqref{P1P2-2}, for each prime $q\in (P_1/2,P_1]$, there is
a multiset $T=S_q$ of residues modulo $q$ such that
\be\label{SqFt}
|f_{S_q}| \le 15 \frac{\log (P_1)}{P_0}:=\eps_1.
\ee
We have $|S_q| =S$ for all $q$, where $S=R_0\sum_{P_0/2<p\le P_0} p$.
Now define a multiset $\{z_1,\dots,z_n\}$ as a union of
multisets $\{e(s/q): s\in S_q, q\in (P_1/2,P_1]\}$. We have,
for $1\le k\le N$,
$$
\sum_{j=1}^n z_j^k=\sum_{P_1/2<q \le P_1} B(q,k), \qquad
B(q,k) = \sum_{s\in S_q} e\pfrac{ks}{q}.
$$
If $q|k$, then $B(q,k)=S$.
When $q \nmid k$, by \eqref{SqFt}, $|B(q,k)| \le \eps_1 S$. Therefore,
\be\label{nodivide}
\sum_{q \nmid k}|B(q,k)| \le \eps_1 n.
\ee
The sum over $q|k$ is estimated at the same way as in Lemma \ref{iter2}:
\be\label{divide}
\sum_{q|k}|B(q,k)| \le \frac{\log N}{\log (P_1/2)} S.
\ee
Combining \eqref{nodivide}, \eqref{divide} and using Proposition \ref{primes}
we arrive at
$$
\frac1n\left|\sum_{j=1}^n z_j^k\right|
\le \eps_1 +  \frac{5\log N}{2P_1},
$$
as required. Moreover, by Proposition \ref{primes} we have
$$
n \le (0.76)^2 R_0 \frac{P_1 P_0^2}{(\log P_0)(\log P_1)}.
$$
Now we take $R_0, P_0, P_1$ the same as in the proof of Theorem \ref{thm1}
so that \eqref{errors2} follows immediately.  The condition \eqref{KN2}
implies \eqref{P1P2-2} for large enough $N$.
\end{proof}

\begin{remark}
As in \cite{5}, one can construct thin sets $T$ modulo $N$ with $|T|=o(L_1L_2)$
and $|f_T|$ small, by iterating Lemma \ref{iter2}.  Roughly speaking, applying
 Lemma \ref{iter4} followed by $r$ iterations of Lemma \ref{iter2}
 produces sets $T$,
with small $|f_T|$, as small as $|T|=O(L_1 L_{r+1})$, where $L_j$ is the $j$-th
iterate of the logarithm of $N$.  We omit the details.
\end{remark}

{\bf Acknowledgments.}  The authors thank Ronald DeVore, Zeev Dvir, 
Venkatesan Guruswami, Piotr Indyk, Sina Jafarpour, Boris Kashin, Howard Karloff,
Imre Leader, Igor Shparlinski and Avi Wigderson for helpful conversations.

%
%
%


\begin{thebibliography}{99}
%


\bibitem{5} M. Ajtai, H. Iwaniec, J. Koml\'os, J. Pintz and
  E. Szemer\'edi, {\it Construction of a thin set with small Fourier
    coefficients}, Bull. London Math. Soc. {\bf 22} (1990), 583--590.

\bibitem{Alon} N. Alon, O. Goldreich, J. H{\aa}stad and R. Peralta, {\it Simple constructions of almost $k$-wise independent random variables},
Random Structures and Algorithms {\bf 3} (3) (1992), 289--303.

\bibitem{An1} J. Andersson, {\it Explicit solutions to certain inf max
problems from Tur\'an power sum theory}, Indag. Math.(N.S.) {\bf 18}
(2007), no. 2, 189--194.

\bibitem{An2} J. Andersson, {\it On the solution to a power sum
problem}, Analytic and probabilistic methods in number theory//
Analiziniai ir tikimybiniai metodai skai\v ci/polhk u teorijoje, 1--5,
TEV, Vilnius, 2007.

\bibitem{An3} J. Andersson, {\it On some power sum problems of Montgomery and
Tur\'an}, Int. Math. Res. Not. IMRN (2008), no. 8, Art. ID rnn015, 9 pp.

\bibitem{BDDW} R. Baraniuk, M. Davenport, R. DeVore and M. Wakin,
{\it A simple proof of the restricted isometry property for random
matrices}, Constr. Approx. {\bf 28}  (2008),  no. 3, 253--263.

\bibitem{BenTa}  A. Ben-Aroya and A. Ta-Shma, {\it Constructing small-bias
sets from algebraic-geometric codes}, Proceedings of the IEEE Symposium on Foundations
of Computer Science (FOCS 2009), IEEE Computer Soc. (2009), 191--197.

\bibitem{BGIKS} R. Berinde, A. Gilbert, P. Indyk, H. Karloff and M. Strauss,
{\it Combining geometry and combinatorics: a unified approach to
  sparse signal recovery}, Proceedigns of the
46th Annual Allerton Conference on Communication, Control and
Computing (2008), 798--805.

\bibitem{BL}  B. Bollobas and I. Leader, {\it Sums in the grid},
Discrete Math. {\bf 162} (1996), 31--48.

\bibitem{Bo} J. Bourgain, {\it Multilinear exponential sums in prime fields
under optimal entropy condition on the sources}, Geom. Funct. Anal. {\bf 18}
(2009), no. 5, 1477-1502.

\bibitem{BG} J. Bourgain and M. Z. Garaev, {\it On a variant of sum-product
estimates and explicit exponential sum bounds in finite fields},
Math. Proc. Cambridge Philos. Soc. {\bf 146} (2009), no. 1, 1--21.

\bibitem{BGlib} J. Bourgain and A. A. Glibichuk, {\it Exponential sum estimate over subgroup in an
arbitrary finite field},  preprint, 2010.

\bibitem{Ca} E. J. Cand\`es, {\it The restricted isometry property and its
implications for compresses sensing}, C. R. Math. Acad. Sci. Paris
{\bf 346} (2008), 589--592.

\bibitem{CRT} E. J. Cand\`es, J. Romberg and T. Tao,
{\it Stable signal recovery from incomplete and inaccurate
measurements}, Comm. Pure Appl. Math. {\bf 59} (2006), 1208--1223.

\bibitem{CT} E. J. Cand\`es and T. Tao, {\it Decoding by linearin Conference in Modern Analysis and Probability 
  programming}, IEEE Trans. Inform. Th. {\bf 51} (2005), 4203--4215.

\bibitem{CT2} E. J. Cand\`es and T. Tao, {\it Near-optimal signal recovery
 from random projections: universal encoding strategies}, IEEE Trans.
 Inform. Theory {\bf 52} (2006), no. 2, 489--509.

\bibitem{D} R. DeVore, {\it Deterministic constructions of compressed
sensing matrices}, Journal of Complexity {\bf 23} (2007), 918--925.

\bibitem{DET} D. Donoho, M. Elad and V. N. Temlyakov,
 {\it On the Lebesgue type inequalities for greedy approximation},
  J. Approximation Theory {\bf 147} (2007), 185--195.

\bibitem{ER} P. Erd\H os and A. R\'enyi, {\it A probabilistic
approach to problems of Diophantine approximation}, Illinois J. Math.
{\bf 1} (1957),  303--315.

\bibitem{GMS} A. C. Gilbert, S. Mutukrishnan and M. J. Strauss,
{\it Approximation of functions over redundant dictionaries using
  coherence}, The 14th Annual ACM-SIAM Symposium on Discrete
Algorithms, (2003), 243--252.

\bibitem{Glu} E. D. Gluskin, {\it An octahedron is poorly approximated by random
subspaces}, Funktsional. Anal. i Prilozhen. {\bf 20} (1986), no. 1, 14--20, 96.

\bibitem{GH} S. Gurevich and R. Hadani, {\it The statistical restricted
 isometry property and the Wigner semicircle distribution of incoherent
 dictionaries}, preprint, arXiv:0812.2602

\bibitem{GS} V. Guruswami and M. Sudan, {\it List decoding algorithms
 for certain concatenated codes}, Proc. 32nd Ann. ACM Sympos. on Theor. 
Computer Sci. (Portland, OR, 2000), 181--190.

\bibitem{IR} K. Ireland, M. Rossen, {\it A classical introduction to
modern number theory}, Springer - Verlag, 1982.


\bibitem{JL} W. B. Johnson and J. Lindenstrauss, {\it Extensions of
Lipschitz mappings into a Hilbert space}, in Conference in Modern Analysis and
Probability (New Haven, Conn., 1982), Contemp. Math. {\bf 26},
Amer. Math. Soc., Providence, 1984, 189--206.

\bibitem{Ka} B. S. Kashin, {\it On widths of octahedron}, Uspekhi Matem.
Nauk {\bf 30} (1975), 251--252 (Russian).

\bibitem{Ka2} B. S. Kashin, {\it Widths of certain finite-dimensional
sets and classes of smooth functions}, Izv. Akad. Nauk SSSR, Ser. Mat. {\bf 41}
(1977), 334--351; English transl. in Math. USSR Izv. {\bf 11} (1978),
317--333.


\bibitem{Kat} N. M. Katz, {\it An estimate for character sums}, J.
Amer. Math. Soc. {\bf 2} (1989), no. 2, 197--200.

\bibitem{Lev} V. I. Levenshtein,
{\it Bounds for packings of metric spaces and some of their applications.} (Russian) Problemy Kibernet. No. 40 (1983), 43–-110.

\bibitem{LT} E. Liu and V. N. Temlyakov, {\it Orthogonal super greedy 
algorithm  and  applications  in  compressed  sensing}, preprint, 2010.

\bibitem{Liv} E. Livshitz, {\it On efficiency of Orthogonal Matching Pursuit},
preprint, 2010, ArXiv: 1004.3946.

\bibitem{MPT} S. Mendelson, A. Pajor and N. Tomczak-Jaegermann,
{\it Reconstruction and subgaussian operators in asymptotic geometric
analysis}, Geom. Funct. Anal. {\bf 17} (2007), 1248--1282.

\bibitem{Mon} H. L. Montgomery, {\it Ten lectures on the interface between
analytic number theory and harmonic analysis}, CBMS, no. 84, 1994.

\bibitem{NR} D. Needle and R. Vershynin, {\it Uniform uncertainty principle
and signal recovery via regularized orthogonal matching pursuit},
Found. Comput. Math. {\bf 9} (2009), no. 3, 317--334.

\bibitem{NT} J. Nelson and V. N. Temlyakov,
{\it On the size of incoherent systems}, preprint, 2010.

\bibitem{Pre} A. Pr\'ekopa, {\it Logarithmic concave measures with application
to statistic processes}, Acta Scient. Math. {\bf 32} (1971), 301--316.

\bibitem{PS}  G. I. Perel'muter and I. E. Shparlinski, {\it Distribution
of primitive roots in finite fields}, Russian Math. Surveys
{\bf 45} (1990),  223--224.

\bibitem{RSW}
 A. Razborov, E. Szemer\'edi and A. Wigderson, {\it
Constructing small sets that are uniform in arithmetic progressions.}
Combin. Probab. Comput. {\bf 2} (1993), no. 4, 513--518.

\bibitem{RS}  J. B. Rosser and L. Schoenfeld, {\it Approximate formulas
for some functions of prime numbers}, Illinois J. Math.
{\bf 6} (1962),  64--94.

\bibitem{RV}  M. Rudelson and R. Vershynin, {\it On Sparse Reconstruction
 from Fourier and Gaussian Measurements}, Comm. Pure Appl. Math. {\bf 61}
 (2008), no. 8, 1025--1045.

\bibitem{St1} S. B. Stechkin, {\it Some extremal problems of trigonometric
sums}, Math. Notes {\bf 55} (1994), no. 1--2, 195--203.

\bibitem{St2} S. B. Stechkin, {\it The Tur\'an problem for trigonometric
sums}, Proc. Steklov Inst. Math. 1994), no. 4 (219), 329--333.

\bibitem{Tao} T. Tao, {\it Open question: deterministic uup matrices},
Weblog at http://terrytao.wordpress.com (2007, July 02).

\bibitem{TaoVu} T. Tao, V. Vu, {\it Additive Combinatorics},
Cambridge University Press, 2006.

\bibitem{Tu} P. Tur\'an, {\it On a new method of analysis and its
applications}, John Wiley \& Sons, Inc., New York, 1984.

\bibitem{Vin} I. M. Vinogradov, {\it An introduction to the theory of numbers},
Pergamon Press, London, New York, 1955.

\bibitem{Wood} D. R. Woodall, {\it A theorem on cubes}, Mathematika
{\bf 24} (1977), 60--62.

\end{thebibliography}
\end{document}